\documentclass[11pt]{article}
\usepackage{amsmath,amsthm,amsfonts,amscd,bezier}
\usepackage{amssymb}
\usepackage{enumerate}
\usepackage[T1]{fontenc}
\usepackage[usenames]{color}

\newtheorem{lemma}{Lemma}[section]

\newtheorem{theorem}[lemma]{Theorem}
\newtheorem{remark}[lemma]{Remark}
\newtheorem{proposition}[lemma]{Proposition}
\newtheorem{corollary}[lemma]{Corollary}
\newtheorem{definition}[lemma]{Definition}

\newcommand{\Dem}{\noindent{\sc Proof:\ \ }}
\newcommand{\cqd}{{\hfill $\rule{2mm}{2mm}$}\vspace{1cm}}

\title{The Analytic Classification of Plane Curves }

\author{Marcelo Escudeiro Hernandes  \\ and \\  Maria Elenice Rodrigues Hernandes\thanks{The authors were partially supported by grant 2019/07316-0, S\~ao Paulo Research Foundation (FAPESP) and the first author was partially supported by CNPq-Brazil Proc. 303638/2020-6.}}

\begin{document}

\maketitle
\markboth{M. E. Hernandes and M. E. Rodrigues
Hernandes}{The Analytic Classification of Plane Curves}

\begin{center} 2010 Mathematics Subject Classification: 14H20 (primary),
14Q05, 14Q20, 32S10 (secondary)\end{center}

\begin{center} keywords: plane curves, analytic classification, analytic invariants. \end{center}

\noindent {}

\begin{abstract}
In this paper, we present a solution to the problem of the analytic
classification of germs of plane curves
with several irreducible components.
 Our algebraic
approach follows precursive
ideas of Oscar Zariski and as a subproduct allows us to recover some particular cases found in the literature.
\end{abstract}

\section{Introduction}

A remarkable step towards the understanding of the local structure of a germ of complex plane curve was the characterization of its topological type. Let $\mathcal{C}_1$ and $\mathcal{C}_2$ be germs of analytic reduced plane curves at the origin of $\mathbb{C}^2$. We say that $\mathcal{C}_1$ and $\mathcal{C}_2$ are topologically equivalent (as embedded germs) if there exist $U$ and $V$ neighborhoods at the origin of $\mathbb{C}^2$ and a homeomorphism $\Phi:U \to V$ such that $\Phi(\mathcal{C}_1 \cap U)= \mathcal{C}_2 \cap V$. In this case, Zariski in \cite{zariski3} and \cite{zariskibook} says that $\mathcal{C}_1$ and $\mathcal{C}_2$ have the same \textit{topological type} or they are \textit{equisingular}. When $\Phi$ is an analytic isomorphism, the curves are called \textit{analytically equivalent}.

The local topology of plane curves has been studied since the first decades of the last century, with important contributions of K. Brauner, W. Burau, O. Zariski, J. Milnor, and others (see, \cite{BriesKnorrer}). For an irreducible plane curve (branch) the topological type of the curve is equivalent to the topology of the complement of the associated link, that is, the intersection of the curve with a small sphere centered at the origin. This intersection is an iterated torus knot. The topological type of the curve is completely described by the type of the knot which is characterized by pairs of integers called Puiseux pairs. There are many complete discrete invariants that determine and are determined by Puiseux pairs, for instance, the characteristic exponents, the value semigroup of the curve, the multiplicity sequence associated to the canonical resolution, and others. For plane curves with several branches, the local topology can be described by the value semigroup of the curve, or equivalently by the value semigroup of each branch and the intersection multiplicities of pairs of branches (see \cite{waldi} and \cite{zariski3}).

With regard to analytic equivalence, there are many hard problems. For instance, how we can decide if two curves are analytically equivalent? What is the associated moduli space?

Any introduction to the analytic classification subject of plane curves is hardly comparable to the review by Washburn presented in \cite{Washburn}. We just cite some contributions to the theme in chronological order:

In 1965, Ebey in \cite{ebey} presented normal forms for some classes of irreducible curves given by parameterizations. In a course given at the Centre de Math\'{e}matiques de l'\'Ecole Polytechnique in 1973, Zariski in \cite{zariskibook} considers  an irreducible analytic plane curve $\mathcal{C}$ and denotes by $\mathbb{L}=\mathbb{L}(\mathcal{C})$ the equisingular class of $\mathcal{C}$, that is, the set of all branches equisingular to $\mathcal{C}$. The moduli space $\mathbb{M}$ is the quotient space of $\mathbb{L}$ by the analytic equivalence relation. Zariski studied the moduli space for some equisingular classes and gave a formula for the dimension of the generic component $\mathbb{M}_g$ of $\mathbb{M}$ for curves with semigroup $\langle n,m\rangle$ in which $1<n < m$ and $m \equiv 1 \ mod \ n$. Few years later, Delorme (in 1978, see \cite{delorme}), considering irreducible curves with semigroup $\langle n,m\rangle$, presented a combinatorial method to obtain the dimension of the generic component of the moduli space.
In 1979, Granger presented normal forms for plane curves with nonsingular transversal branches (an ordinary multiple point) and gave a formula to compute the dimension of $\mathbb{M}_g$ (see \cite{Granger}).

Bruce and Gaffney in 1982 (see \cite{bruce}), classified the simple irreducible plane curves, that is, when the moduli space is a zero dimensional space. A more general approach for the moduli space was presented by Laudal and Pfister in 1988 (see \cite{laudal}) where, for irreducible plane curves with semigroup $\langle n,m\rangle$, they fixed an analytic invariant (the Tjurina number) in order to describe normal forms given by elements in $\mathbb{C}[X,Y]$. Greuel and Pfister (see \cite{Greuel}, 1994) developed a general method to construct coarse moduli spaces for singularities in the sense of Mumford Geometric Invariant Theory. Kang (see \cite{Kang}, 2000), C\^amara and Sc\'ardua (see \cite{leonardo}, 2018), using different methods, classified analytic plane curves defined by weighted homogeneous polynomials. Kolgushkin and Sadykov in 2001 (see \cite{KolgSad}) obtained normal forms for stably simple reducible curve singularities in complex spaces of any dimension.

In 2011, Hefez and the first author (see \cite{HefezHern}) presented a solution for the analytic classification of irreducible plane curves as considered by Zariski in \cite{zariskibook} (see also \cite{handbook}). Stratifying each equisingularity class by the analytical invariant given by the set of values of K\"ahler differentials, they
presented normal forms for each stratum that allowed them to describe the moduli space. Later, the authors in a joint work with Hefez (in 2015, see \cite{HefezHernRod}), generalized such method to obtain the analytic classification of plane curves with two irreducible components.

Genzmer and Paul (in 2016, see \cite{genzmerpaul}), using tools of Foliation theory, described the moduli space for generic plane curves such that every branch admits semigroup $\langle n,m\rangle$ and they presented a method to obtain the normal form for the generic case. Ayuso and Rib\'on (in 2020, see \cite{pedro}), using holomorphic flows recovered the normal forms for branches presented in \cite{HefezHern}. In 2022, Genzmer in \cite{Genzmer}, obtained a formula for the dimension of $\mathbb{M}_g$ for any irreducible plane curve.

It is remarkable to observe that for about fifty years the questions proposed by Oscar Zariski in \cite{zariskibook}, directly or indirectly, motivated these and so many other researchers, that we could not mention. The aim of this work is to present an answer to the intricate problem
of the analytic classification of reduced plane curves in a fixed topological class. The combinatorial issues involved in this problem were probably our biggest challenge and  makes this case not a simple generalization of the results about curves with two irreducible components presented in \cite{HefezHernRod}.

We consider a topological class fixing the \textit{value semiring} $\Gamma$ of a reduced plane curve $\mathcal{C}$. The value semiring is equivalent, as a set, to the classical value semigroup. The main difference is that the value semiring, equipped with the tropical operations, admits a finite minimal set of generators in contrast to the value semigroup (see \cite{carvalho1}), and such generators allow us to recover directly the topological data pointed out by Zariski, that is, the value semigroup of each branch and their mutual intersection multiplicities. With our approach we recover the normal forms for the irreducible and two branches cases, presented in \cite{HefezHern} and \cite{HefezHernRod}, but it is not a simple induction step. In fact, the passage $1 \to 2 \to r$ branches requires extra efforts and finer analysis in several situations. For example, understanding the semigroup for two branches as presented by Garcia (see \cite{Garcia}) is not sufficient to describe the semigroup for $r\geq 3$. In fact, Mata (see \cite{Mata}) introduces new ingredients and non trivial combinatorial aspects to perform this task. This important object and other invariants defined by a set of values of fractional ideals are described in Section 2.

Our strategy is, under the action of a permutation group, to consider the branches of the curve in blocks according to their tangent lines, then we order the blocks and their elements by the multiplicities of the branches. In addition, we identify the appropriated group $\mathcal{G}$ of local diffeomorphisms that preserve such properties. This will be done in Section \ref{section-groupactions}. In Section 4 we obtain the Puiseux block form of a multigerm (Proposition \ref{prop-puiseux-block}) that is a convenient parameterization of each component of the multigerm.

Section 5 contains the main results of this work. We introduce the $\mathcal{G}$-invariant $\Lambda_{\mathcal{G}}$ that corresponds to values of elements in a fractional ideal $\mathcal{I}_{\mathcal{G}}$ of the local ring $\mathcal{O}$ of $\mathcal{C}$. The elements of $\mathcal{I}_{\mathcal{G}}$ are closely related to elements of the tangent space to the $\mathcal{G}$-orbit of the multigerm $\varphi$ associated to the curve $\mathcal{C}$ (Proposition \ref{DifKTg}). In Theorem \ref{A1normal-form} we present a {\it $\mathcal{G}$-normal form} of $\varphi$ by a reduction process using the set $\Lambda_{\mathcal{G}}$, more specifically by their fibers. Theorem \ref{Anormal-form} takes into account the homothety group action and by this result we can decide if two curves are analytically equivalent or not as we discuss in the end of the section.

Finally, in Section 6 we apply our techniques to recover some known
results related with the analytic equivalence of plane curves: the
irreducible case, bigerms and the ordinary multiple point
singularity. Concerning the case of plane curves in which each branch
admits value semigroup $\langle n,m\rangle$ and mutual intersection
multiplicities equal to $nm$ we apply our results to obtain a
pre-normal form taking into account the value semiring $\Gamma$
(Proposition \ref{pre-normal}) and we discuss about the generic
component of the corresponding moduli space (Proposition
\ref{generic-dim} and Corollary \ref{23}).

\section{Analytic equivalence and invariants}
\label{section:AnalInvariants}

We denote by $\mathbb{C}\{X,Y\}$ the power series ring with complex coefficients in the variables $X$ and $Y$, which are absolutely convergent in a neighborhood of the origin in $\mathbb{C}^2$ and by $\mathcal{M}=\langle X,Y \rangle$ its maximal ideal. Let $\mathcal{C}$ be the germ of an analytic reduced plane curve in $(\mathbb{C}^2,0)$ defined by a reduced element $f\in \mathcal{M}\setminus\{0\}$, or in other words, $\mathcal{C}=f^{-1}(0)\cap U$, where $U$ is an open neighborhood at the origin $0$ in $\mathbb{C}^2$.

If $f=f_1 \cdot \ldots \cdot f_r$ is the decomposition of $f$ into irreducible factors, then each $f_i$ defines a branch, that is, an irreducible analytic plane curve, denoted by $\mathcal{C}_i$ for $i=1, \ldots, r$.

Given a branch $\mathcal{C}_i$ of $\mathcal{C}$ with multiplicity $n_i$, that is, $f_i\in \mathcal{M}^{n_i}\setminus\mathcal{M}^{n_i+1}$, we can parameterize it by $(a_it_i^{n_i}+\ldots, b_it_i^{n_i}+\ldots)\in \mathbb{C}\{t_i\}\times\mathbb{C}\{t_i\}$ with $a_i\neq 0$ or $b_i\neq 0$. In what follows we will consider primitive parameterizations that is those that cannot be obtained from other ones composed with higher powers of $t_i$. We will often identify such parameterization with the smooth map-germ $\varphi_i:(\mathbb{C},0)\to
(\mathbb{C}^2,0)$ defined by $t_i\mapsto (x_i,y_i):=(a_it_i^{n_i}+\ldots, b_it_i^{n_i}+\ldots)$.

We call
$$[\varphi_1, \dots ,\varphi_r]=\left[\begin{array}{ccccc}
x_1  & \cdots & x_i  & \cdots & x_r \\
y_1  & \cdots  & y_i & \cdots & y_r \\
\end{array}\right]$$ a multigerm for $\mathcal{C}$ and we denote by $\mathcal{P}$ the set of multigerms of plane curves with $r$ branches.

The analytic equivalence of reduced plane curves with $r$ branches is translated, as we noticed in \cite{HefezHernRod}, into $\mathcal{S}\times{\cal A}$-equivalence on $\mathcal{P}$, where $\mathcal{S}$ is
the symmetric group on $r$ elements and ${\cal A}=\{(\rho_1, \ldots, \rho_r , \sigma);  \ \rho_i \in
{\rm Diff}(\mathbb{C},0) \ \mbox{and} \ \sigma \in {\rm Diff}(\mathbb{C}^2,0), \ 1 \leq i \leq
r\}$ is the group of right-left equivalence, where ${\rm Diff}(\mathbb{C}^l,0)$ denotes the diffeomorphism group of $(\mathbb{C}^l,0)$.

Given $(\rho_1, \ldots, \rho_r, \sigma) \in \mathcal{A}$ and  $\varphi=[\varphi_1, \dots ,\varphi_r] \in \mathcal{P}$ the action of $\mathcal{A}$ on $\mathcal{P}$ is defined as:
{\small $$(\rho_1, \ldots, \rho_r, \sigma)\cdot \varphi:= [\sigma \circ \varphi_1 \circ \rho_1^{-1}, \ldots,
 \sigma \circ \varphi_r \circ \rho_r^{-1}]=\left[\begin{array}{ccc}
\sigma_1 (x_1 \circ \rho_1^{-1}, y_1 \circ \rho_1^{-1})  & \cdots  & \sigma_1 (x_r \circ \rho_r^{-1}, y_r \circ \rho_r^{-1}) \\
\sigma_2 (x_1 \circ \rho_1^{-1}, y_1 \circ \rho_1^{-1})  & \cdots  & \sigma_2 (x_r \circ \rho_r^{-1}, y_r \circ \rho_r^{-1}) \\
\end{array}\right],$$
}
where $\sigma=(\sigma_1,\sigma_2)$.

We say that
 $\varphi,\psi \in \mathcal{P}$ are $\mathcal{A}$-equivalent, denoted by $\varphi \underset{\mathcal{A}}{\sim}\psi$, if and only if they are in the same $\mathcal{A}$-orbit. Now, the group $\mathcal{S}\times\mathcal{A}$ acts on $\mathcal{P}$ by
 $$\left( \pi, (\rho_1, \ldots, \rho_r, \sigma)\right)\cdot \varphi:= \pi ((\rho_1, \ldots, \rho_r, \sigma)\cdot \varphi):=[\sigma \circ \varphi_{\pi(1)} \circ \rho_{\pi(1)}^{-1}, \ldots,
 \sigma \circ \varphi_{\pi(r)} \circ \rho_{\pi(r)}^{-1}],$$
 for all $\pi \in \mathcal{S}, \ (\rho_1, \ldots, \rho_r, \sigma) \in \mathcal{A}$ and $\varphi \in \mathcal{P}$.

It is a well known fact that for a plane curve $\mathcal{C}$ defined by $f\in\mathbb{C}\{X,Y\}$, the isomorphism class of its local ring $\mathcal{O}:=\frac{\mathbb{C}\{X,Y\}}{\langle f\rangle}$ completely determines the analytic class of $\mathcal{C}$. The aim of this section is to present discrete $\mathcal{A}$-invariants related to some fractional ideals of $\mathcal{O}$.

A parameterization $\varphi_i =(x_i,y_i)\in\mathbb{C}\{t_i\}\times\mathbb{C}\{t_i\}$ of a branch $\mathcal{C}_i$ given by $f_i\in\mathbb{C}\{X,Y\}$ provides the exact sequence
$$\{0\}\rightarrow \langle f_i\rangle\rightarrow \mathbb{C}\{X,Y\}\underset{\varphi_i^*}{\rightarrow} \mathbb{C}\{x_i,y_i\}\rightarrow \{0\},$$ where  $\varphi_i^*(h):=h(x_i,y_i)$.
In what follows, we will identify
$\mathcal{O}_i:=\frac{\mathbb{C}\{X,Y\}}{\langle f_i\rangle}$ with the subalgebra
$\mathbb{C}\{x_i,y_i\}\subseteq \mathbb{C}\{t_i\}$. In particular, the integral closure of
$\mathcal{O}_i$ in its quotient field $\mathcal{Q}_i=\mathbb{C}((t_i))$ is
$\overline{\mathcal{O}_i}=\mathbb{C}\{t_i\}$ and
$\mathcal{O}_i=\overline{\mathcal{O}}_i$ if and only if $n_i=1$.

We consider the natural discrete normalized valuation
$$\begin{array}{cccl}
\nu_i: & \mathcal{Q}_i & \longrightarrow & \overline{\mathbb{Z}}:=\mathbb{Z}\cup\{\infty\} \\
& \frac{p}{q} & \longmapsto & \nu_i\left (\frac{p}{q}\right ):=ord_{t_i}(p)-ord_{t_i}(q),
\end{array}$$
where $p,q\in \mathcal{O}_i, q\neq 0$ and
$\nu_i(0)=\infty$. The set $$\Gamma_i=\{\nu_i(p);\
p\in\mathcal{O}_i\}=\{\nu_i(h):=\nu_i(\varphi^*_i(h));\
h\in\mathbb{C}\{X,Y\} \}$$ is a submonoid of
$\overline{\mathbb{N}}:=\mathbb{N}\cup\{\infty\}$ and
$S_i=\Gamma_i\cap\mathbb{N}$ is the classical {\it value semigroup} of the branch $\mathcal{C}_i$, which is an $\mathcal{A}$-invariant
and a complete topological invariant.

The conductor ideal of $\mathcal{O}_i$ in $\overline{\mathcal{O}_i}$ is
$(\mathcal{O}_i:\overline{\mathcal{O}_i})=\{h \in \mathcal{O}_i;\ h \overline{\mathcal{O}_i}\subseteq \mathcal{O}_i\}=\langle
t_i^{\mu_i}\rangle$, where $\mu_i$, called conductor of $\Gamma_i$, satisfies $\mu_i-1\not\in\Gamma_i$ and
$\mu_i+\overline{\mathbb{N}}\subset\Gamma_i$. In this case, the conductor of $\Gamma_i$ coincides with the Milnor
number of $\mathcal{C}_i$, that is, $\mu_i=\dim_{\mathbb{C}}\frac{\mathbb{C}\{X,Y\}}{\langle
    (f_i)_X,(f_i)_Y\rangle}$ where $(f_i)_X$ and $(f_i)_Y$ indicate the derivatives of $f_i$ with respect to $X$ and $Y$ respectively.

The previous concepts can be extended for a reduced curve
$\mathcal{C}$ with $r$ branches defined by $f=f_1\cdot\ldots\cdot
f_r$.

Considering $I=\{1,\ldots ,r\}$ and the monomorphism
$$\begin{array}{ccc}
\mathcal{O} & \rightarrow & \bigoplus_{i\in I}\mathcal{O}_i \\
p & \mapsto & (p_1,\ldots ,p_r),\end{array}$$ where $p_i$ denotes the class of $p\in \mathcal{O}$ in $\mathcal{O}_i$, it is possible to verify that the total ring of fractions of $\mathcal{O}$ is $\mathcal{Q}=\bigoplus_{i\in I}\mathcal{Q}_i$ and the integral closure of $\mathcal{O}$ in $\mathcal{Q}$ is $\overline{\mathcal{O}}=\bigoplus_{i\in I}\mathbb{C}\{t_i\}$. As before, if
$\varphi=[\varphi_1,\ldots, \varphi_r]\in\mathcal{P}$ is a multigerm associated to
$\mathcal{C}$ then we identify
$\mathcal{O}$ with the subalgebra $\{\varphi^*(h):=(\varphi^*_1(h),\ldots
,\varphi^*_r(h)),\ h\in \mathbb{C}\{X,Y\}\}\subseteq \bigoplus_{i\in I}\mathcal{O}_i$ and we set
 $$\Gamma = \left\{\nu(h):=\left(\nu_1(h),\dots,\nu_r(h)\right); \ h \in \mathbb{C}\{X,Y\} \right\}\subset \Gamma_1\times\ldots \times\Gamma_r.$$
The conductor ideal of $\mathcal{O}$ in $\overline{\mathcal{O}}$ is $$(\mathcal{O}:\overline{\mathcal{O}})=\bigoplus_{i\in I}\left ((\mathcal{O}_i:\overline{\mathcal{O}_i})\cdot  \prod_{{j\in I}\atop {j\neq i}}\varphi_i^*(f_j) \right )=\bigoplus_{i\in I}\langle t_i^{\kappa_i}\rangle,$$ where
$\kappa_i=\mu_i+ \sum_{{j\in I\atop j \neq i} }\nu_i(f_j)$. The element $\kappa=(\kappa_1, \ldots, \kappa_r)$, called conductor of $\Gamma$, satisfies $\kappa + \overline{\mathbb{N}}^r \subset \Gamma$ and $\kappa-e_i\not\in\Gamma$ for any element $e_i$ in the canonical basis of $\mathbb{Q}^r$. For a reference of these facts we indicate \cite{Gorenstein}.

Zariski in \cite{zariski3} showed that (up to permutation of the branches) the topological class of $\mathcal{C}$ is totally
characterized by $\Gamma_i\cap \mathbb{N}$ and the intersection multiplicity $\nu_i(f_j)\ (=\nu_j(f_i))$ of $\mathcal{C}_i$ and $\mathcal{C}_j$ for
$i,j\in I$ with $i\neq j$. On the other hand, Waldi (see
\cite{waldi}) proved that $\Gamma\cap\mathbb{N}^r$ is also a complete
topological invariant and Mata, in \cite{Mata},
showed that $\Gamma$ determines and it is determined by $\Gamma_i$
and $\nu_i(f_j)$ for all $i,j\in I,\ i\neq j$, connecting the results of Waldi and Zariski.

In contrast to $\Gamma_i$, the (additive) semigroup $\Gamma$ is not
finitely generated, but equipped with the tropical operations, the set $(\Gamma,\inf,+)$ is a
finitely generated semiring\footnote{The set $(\Gamma,\inf,+)$ is a semiring since $(\Gamma,\inf)$ and  $(\Gamma,+)$ are monoids with identity elements $\underline{\infty}=(\infty,\ldots ,\infty)$ and $\underline{0}=(0,\ldots ,0)$ respectively; $\inf\{\alpha+\beta,\alpha+\gamma\}=\alpha+\inf\{\beta,\gamma\}$ and $\underline{\infty}+\alpha=\underline{\infty}$ for every $\alpha,\beta,\gamma\in\Gamma$.} (see \cite{carvalho1}), where $\inf\{(\alpha_1,\ldots ,\alpha_r),(\beta_1,\ldots
,\beta_r)\}=(\min\{\alpha_1,\beta_1\},\ldots
,\min\{\alpha_r,\beta_r\})$, that we call the {\it value semiring} associated to $\mathcal{C}$. More precisely,
there exists a minimal set of generators $\{v_j,\ j=1,\ldots
,g\}\subset\Gamma$ such that any $\gamma\in\Gamma$ can be written as
\begin{equation}\label{semiring}\gamma=\inf\left \{\sum_{j=1}^{g}a_{1j}v_j,\ldots,\sum_{j=1}^{g}a_{rj}v_j\right \},\end{equation}
where $a_{kj}\in\mathbb{N},\ k\in I,\ 1\leq j\leq g$.
The elements in $\{v_j,\ j=1,\ldots ,g\}\cap\mathbb{N}^r$ are
precisely the valuation of branches that achieve maximal contact with
some $\mathcal{C}_i$ for $i\in I$ and the subset of generators with some
coordinate equal to $\infty$ is precisely $\{\nu(f_i),\ i\in I\}$
(see \cite{Mata} and \cite{carvalho1}).

In a more general situation, given any (regular) fractional ideal $\mathcal{I}\subseteq \mathcal{Q}=\bigoplus_{i\in I}\mathcal{Q}_i$ of $\mathcal{O}$ we can consider the set of values $$\nu(\mathcal{I})=\left \{\nu(z):=(\nu_1(z_1),\ldots ,\nu_r(z_r));\ z=(z_1,\ldots ,z_r)\in\mathcal{I}\subseteq\mathcal{Q}\right \}\subseteq \overline{\mathbb{Z}}^r.$$

The set $(\nu(\mathcal{I}),\inf)$ is a $\Gamma$-semimodule, that is,
$\Gamma + \nu(\mathcal{I}) \subseteq \nu(\mathcal{I})$. Moreover, there exists (unique) $(\zeta_1,\ldots ,\zeta_r)=:\inf(\nu(\mathcal{I}))\in\nu(\mathcal{I})$ such that $\zeta_i\leq\alpha_i$ for every $(\alpha_1,\ldots ,\alpha_r)\in\nu(\mathcal{I})$ and $\nu(\mathcal{I})$ admits
a ``conductor'' $\varrho \in \nu(\mathcal{I})$, that is, $\varrho +
\overline{\mathbb{N}}^r \subseteq \nu(\mathcal{I})$ and $\varrho
-e_i\not\in \nu(\mathcal{I})$ for all $e_i$ in the canonical
$\mathbb{Q}$-basis of $\mathbb{Q}^r$ (see \cite{edson} for more properties concerning values set of fractional ideals of $\mathcal{O}$).

The next definition generalizes the corresponding concept introduced
by Garcia in \cite{Garcia} and Mata in \cite{Mata} for
$\Gamma\cap\mathbb{N}^r$.

\begin{definition} Let $\mathcal{I}$ be a fractional ideal of $\mathcal{O}$ with set of
    values $\Delta=\nu(\mathcal{I})$. Given $J\subseteq I$ the $J$-{\bf fiber} of $\gamma=(\gamma_1,\ldots
    ,\gamma_r)\in\overline{\mathbb{N}}^r$ with respect to $\Delta$ is the
    set
    $$F^{\Delta}_{J}(\gamma)=\{(\delta_1,\ldots ,\delta_r)\in\Delta;\ \delta_j=\gamma_j\ \mbox{for every}\ j\in J\ \mbox{and}\ \delta_i>\gamma_i\ \mbox{for all}\ i\in I\setminus J\}.$$

    We say that $\gamma \in\overline{\mathbb{N}}^r$ is
    \begin{enumerate}
        \item[(i)] a \textbf{maximal} element of $\Delta$, if $F^{\Delta}_I(\gamma)\neq \emptyset$ (i.e. $\gamma\in\Delta$) and
        $F^{\Delta}_{\{i\}}(\gamma)= \emptyset$ for all $i\in I$;
        \item[(ii)] a \textbf{relative maximal} element of $\Delta$, if $\gamma$ is maximal and $F^{\Delta}_J(\gamma) \neq \emptyset, \ \ \forall \ J \subseteq I$ with $\sharp J \geq 2$;
        \item[(iii)] an \textbf{absolute maximal} element of $\Delta$, if $\gamma$ is maximal and $F^{\Delta}_J(\gamma) = \emptyset, \ \forall \ J \subsetneq I, \ J \neq
        \emptyset$.
    \end{enumerate}
\end{definition}

\begin{remark}
	\label{Remark-Maximal-Box}
	Let $\mathcal{I}$ be a (regular) fractional ideal of $\mathcal{O}$ such that $\inf(\nu(\mathcal{I}))=(\zeta_1, \ldots, \zeta_r)$ and $\varrho=(\varrho_1, \ldots, \varrho_r)$ is the conductor of $\nu(\mathcal{I})$. In \cite{carvalho2} (Theorem 16) the authors presented an algorithm to compute a finite set of generators of $\nu(\mathcal{I})$ (as $\Gamma$-semimodule) and showed that $\nu(\mathcal{I})$ is characterized by its elements in the box $R=[\zeta_1,\varrho_1]\times \ldots \times [\zeta_r,\varrho_r]$. In particular, we have that all maximal elements of $\nu(\mathcal{I})$ belong to $R$.

For $\mathcal{I}=\mathcal{O}$ we obtain $\nu(\mathcal{O})=\Gamma$ and $\Gamma$ is determined by its elements in $[0,\kappa_1]\times\ldots\times [0,\kappa_r]$.
\end{remark}

Other $\mathcal{A}$-invariants that play an important role in this
work are related with set of orders of differential 1-forms.

Let $\Omega^1=\mathbb{C}\{X, Y\}dX+\mathbb{C}\{X, Y\}dY$ be the $\mathbb{C}\{X,Y\}$-module of differential 1-forms. If $\varphi_i=(x_i,y_i)\in\mathbb{C}\{t_i\}\times\mathbb{C}\{t_i\}$ is a parameterization of a plane branch $\mathcal{C}_i$ defined by $f_i$, we extend the homomorphism $\varphi^*_i$ to $\Omega^1$ in the following way:

Given $\omega=a(X,Y) dX+b(X,Y) dY \in \Omega^1$ we define
$$\varphi^*_i(\omega):=t_i\cdot (\varphi^*_i(a)\cdot x'_i+\varphi^*_i(b)\cdot y'_i)\in\mathcal{Q}_i$$
where $x'_i$ and $y'_i$ denote, respectively, the derivative of $x_i,y_i\in\mathbb{C}\{t_i\}$ with respect to $t_i$.

Notice that the kernel of $\varphi_i^*$ is $$\{\omega\in\Omega^1;\ \varphi^*_i(\omega)=0\}=\left\{\omega\in\Omega^1;\ \frac{\omega\wedge df_i}{dX\wedge dY}\in \langle f_i\rangle\right \}=f_i\cdot\Omega(log\ \mathcal{C}_i)$$
where $\Omega(log\ \mathcal{C}_i)$ is the module of logarithmic differential forms along $\mathcal{C}_i$ which is the dual module of the logarithmic vector field $Der(-log\ \mathcal{C}_i)$ along $\mathcal{C}_i$ (see \cite{Saito}).

In addition, if $\mathcal{F}_i:=f_i\cdot\Omega^1+\mathbb{C}\{X,Y\}\cdot df_i$, then the K\"ahler differential module of $\mathcal{O}_i$ is $\Omega_{i}\approx\frac{\Omega^1}{\mathcal{F}_i}$ and its torsion submodule is given by $\mathcal{T}_{i}\approx\frac{f_i\cdot\Omega(log\ \mathcal{C}_i)}{\mathcal{F}_i}$. Moreover $\frac{\Omega_i}{\mathcal{T}_{i}}\approx\frac{\Omega^1}{f_i\cdot\Omega(log\ \mathcal{C}_i)}\approx \varphi^*_i(\Omega^1)$.

The set $\varphi^*_i(\Omega^1)\subset \mathcal{Q}_i$ is a fractional
ideal of $\mathcal{O}_i$ and
$$\Lambda_i=\{\nu_i(\omega):=\nu_i(\varphi^*_i(\omega));\
\omega\in\Omega^1\}$$ is an $\mathcal{A}$-invariant of $\mathcal{C}_i$.
As $\nu_i(dh)=\nu_i(h)$ for any $h\in\mathcal{M}$ we have
$\Gamma_i\setminus\{0\}\subseteq\Lambda_i$. The
set $\Lambda_i$ is one of the main tools considered in \cite{HefezHern} in order to classify plane
branches up to analytical equivalence.

Similarly for a multigerm  $\varphi=[\varphi_1,\ldots ,\varphi_r]$ we get the $\mathcal{A}$-invariant
\begin{equation}
\label{LambdaInvariant}
\Lambda =\nu(\varphi^*(\Omega^1))=\{\nu(\omega):=(\nu_1(\omega),\ldots ,\nu_r(\omega));\ \omega\in\Omega^1\}\subset\Lambda_1\times\ldots\times\Lambda_r \subset \overline{\mathbb{N}}^r.
\end{equation}
Remark that $\varphi^*(\Omega^1)\approx\frac{\Omega^1}{f\cdot\Omega(log\ \mathcal{C})}$ is a fractional ideal of $\mathcal{O}$ with $\Gamma\setminus\{\underline{0}\}\subseteq\Lambda$ and conductor $\varrho=(\varrho_1,\ldots, \varrho_r)$ satisfying $\varrho_i\leq\kappa_i$ for $i=1,\ldots ,r$.

In particular, by Remark \ref{Remark-Maximal-Box}, the set $\Lambda$ is characterized by its points in the box
$[0,\varrho_1]\times\ldots\times [0,\varrho_r]\subseteq [0,\kappa_1]\times\ldots\times [0,\kappa_r]$. Consequently for each
topological class determined by the semiring $\Gamma$ there is a
finite number of possible $\Lambda$ sets.  

The set $\Lambda$ is related to the Tjurina number of $\mathcal{C}$ (see \cite{tjurina}), the set $\nu(J(f))=\{\nu(h);\ h\in
J(f):=\langle f,f_X,f_Y\rangle\subseteq\mathbb{C}\{X,Y\}\}$ and
the values of residues of elements in $f\cdot\Omega(log\
\mathcal{C})$. More explicitly, if $\omega\in f\cdot\Omega(log\
\mathcal{C})$, then there exist $\eta\in\Omega^1$,
$g,h\in\mathbb{C}\{X,Y\},\ g\not\in\bigcup_{i=1}^{r}\langle
f_i\rangle$ such that $g\cdot\omega=h\cdot df+f\cdot\eta$. The class
$res(\omega)$ of $\frac{h}{g}$ in $\mathcal{Q}$ is called the residue
of $\omega$ and we put $Res(f)=\{res(\omega),\ \omega\in
f\cdot\Omega(log\ \mathcal{C})\}$ (see \cite{Saito}). In \cite{Pol}, Pol showed that
$\nu(J(f))=\Lambda+\kappa -(1,\ldots ,1)$ and $\lambda\in\Lambda$ if
and only if $-\lambda\not\in\nu(Res(f))$.

Using the value semiring $\Gamma$ or more specifically the subset $\nu(\mathcal{M}^2)$ of values of $\mathcal{M}^2$ we can obtain a distinguished presentation for a multigerm as we show in Section $4$. We will consider particular fractional ideals of $\mathcal{O}$ determined by differential $1$-forms that are related to the tangent space to the orbit of a multigerm according to the action of a Lie group.

\section{Group actions}\label{section-groupactions}

Let $\varphi=[\varphi_1, \ldots, \varphi_r]\in\mathcal{P}$ such that $\varphi_i(t_i)=(a_it_i^{n_i}+\ldots, b_it_i^{n_i}+\ldots)$ with $a_i \neq 0$ or $b_i \neq 0$. Given $(\rho_1, \ldots, \rho_r, \sigma)\in\mathcal{A}$ with $\rho_i^{-1}(t_i)=u_it_i+\ldots$ and $\sigma(X,Y) =(\alpha X + \beta Y+\ldots, \gamma X + \delta Y+\ldots)$ such that $u_i, \alpha \delta - \beta \gamma\in\mathbb{C}^*:=\mathbb{C}\setminus\{0\}$ we have
$$\psi_i(t_i):=\sigma \circ \varphi_i\circ\rho_i^{-1}(t_i) =((\alpha a_i+\beta b_i)u_i^{n_i}t_i^{n_i}+\ldots, (\gamma a_i+\delta b_i)u_i^{n_i}t_i^{n_i}+\ldots).$$

Denoting $j^k\xi$ the $k$-jet of a map-germ $\xi$ at the origin we put $j^k\varphi:=[j^k\varphi_1, \ldots, j^k\varphi_r]$.

The invariance of the multiplicity of the branches gives us a one-to-one correspondence between $j^1 \sigma$ and a M\"{o}bius transformation on $\overline{\mathbb{C}}:=\mathbb{C} \cup \{\infty\}$ given by $T(\theta)=\frac{\gamma +\delta \theta}{\alpha +\beta \theta}$ such that if $\theta_i:=\frac{b_i}{a_i}\in \overline{\mathbb{C}}$ is the slope of the tangent line $a_iY-b_iX=0$ of $\mathcal{C}_i$ at the origin, then $T(\theta_i)$ is the slope of the tangent line to the curve corresponding to $\psi_i$ for $1\leq i\leq r$. In this way, up to change of coordinates, we can consider the tangent lines of the branches in a such way that up to three chosen directions are fixed.

If all branches have the
same tangent line at the origin with slope $\theta$, considering $\gamma=-\delta \theta$, then $T(\theta)=0$. If the branches have exactly two tangent lines with distinct slopes $\theta_1$ and $\theta_2$, then taking
$\gamma=-\delta \theta_1$ and $\alpha=-\beta \theta_2$ we obtain $T(\theta_1)=0$ and $T(\theta_2)=\infty$. On the other hand if there exist at least three branches with distinct tangent lines of slopes $\theta_1,\theta_2$ and $\theta_3$, we can consider the unique M\"{o}bius transformation $T$ such that $T(\theta_1)=0$, $T(\theta_2)=\infty$ and $T(\theta_3)=1$.

By the above description and an $\mathcal{S}\times\mathcal{A}$-action, we may adjust the components of $\varphi\in\mathcal{P}$ in a particular way:

\begin{definition}\label{blockform} Given $\varphi=[\varphi_1, \ldots, \varphi_r]\in \mathcal{P}$ we say that $\varphi$ is in {\bf block form} if there are sets
$B_1=\{\varphi_{k_1}=\varphi_1, \ldots, \varphi_{k_2-1}\}, B_2=\{\varphi_{k_2}, \ldots, \varphi_{k_3-1}\},
\ldots ,B_{s}=\{\varphi_{k_{s}}, \ldots, \varphi_{r}\}$
satisfying:
\begin{enumerate}
\item[i)] The tangent line of every branch in $B_i$ has the same slope $\theta_i$ and, if $i\neq j$ then $\theta_i\neq \theta_j$ for $1\leq i,j\leq s$. In addition, $\theta_1=0, \theta_2=\infty$ and $\theta_3=1$.
\item[ii)] If $\varphi_j$ and $\varphi_k$ are elements of $B_i$ with $j<k$, then $n_j\leq n_k$.
\item[iii)] $n_{k_1}\leq n_{k_2}\leq \ldots \leq n_{k_{s}}$.
\end{enumerate}
A set $B_i$ as above is called a {\bf block} of $\varphi$.
\end{definition}

It is clear that any element in $\mathcal{P}$ is $\mathcal{S}\times \mathcal{A}$-equivalent to a multigerm in block form.

From now on we consider $\mathcal{A}$-action on block form multigerms in $\mathcal{P}$ and without loss of generality we can suppose $j^{n_i}\varphi_i=(a_it_i^{n_i}, b_it_i^{n_i})$ such that
$$\begin{array}{lll}
a_i=1, & b_i=0 & \mbox{if}\ \varphi_i\in B_1, \\
a_i=0, & b_i=1 & \mbox{if}\ \varphi_i\in B_2, \\
a_i=1, & b_i=1 & \mbox{if}\ \varphi_i\in B_3\ \mbox{and}\\
a_i=1, & b_i=\theta_j & \mbox{if}\ \varphi_i\in B_j\ \mbox{for}\ j>3.
\end{array}
$$

The next lemma is an immediate consequence of the above explanation and gives us the subgroups of $\mathcal{A}$ that preserve the tangent cone of the curve $\mathcal{C}$ and the properties of the block form.

In what follows we consider the subgroups of $\mathcal{A}$:
$$\begin{array}{l}
\mathcal{H}=\{(\rho_1,\ldots, \rho_r,\sigma)\in \mathcal{A}; \ \rho_i=u_it_i,\ \sigma =(\alpha X, \delta Y)\ \mbox{and}\ u_i,\alpha,\delta \in \mathbb{C}^{*}\},\\
\mathcal{H}'=\{(\rho_1,\ldots, \rho_r,\sigma)\in \mathcal{A}; \ \rho_i=u_it_i,\ \sigma =(\alpha X, \alpha Y)\ \mbox{and}\ u_i,\alpha \in \mathbb{C}^{*}\},\\
\tilde{\mathcal{A}}_1=\{(\rho_1,\ldots, \rho_r,\sigma)\in \mathcal{A};\ j^1\rho_i=t_i,\ j^1\sigma =(X+\beta Y, Y)\ \mbox{and}\ \beta \in \mathbb{C}\}, \\
\mathcal{A}_1=\{(\rho_1,\ldots, \rho_r,\sigma)\in \mathcal{A};\ j^1\rho_i=t_i\ \mbox{and}\  j^1\sigma =(X, Y)\}.\end{array}$$

\begin{lemma}\label{subgroups} The subgroups of $\mathcal{A}$ that preserve the  elements of $\mathcal{P}$ as in Definition \ref{blockform}, according to number $s$ of blocks, are:
\begin{center}\begin{tabular}{|c|c|}
  \hline
 $s$ & Subgroup of $\mathcal{A}$ \\ \hline
  $1$ & $\mathcal{H}\circ \tilde{\mathcal{A}}_1$  \\ \hline
  $2$ & $\mathcal{H} \circ \mathcal{A}_1$  \\ \hline
 $\geq 3$ &  $\mathcal{H}' \circ \mathcal{A}_1$ \\ \hline
\end{tabular}
\end{center}
\end{lemma}
\Dem The description follows from the correspondence between the M\"{o}bius transformation $T(\theta)=\frac{\gamma +\delta \theta}{\alpha + \beta \theta}$ and $j^1 \sigma=(\alpha X + \beta Y, \gamma X + \delta Y)$ with $\sigma \in {\rm Diff}(\mathbb{C}^2,0)$.

For $s=1$, all branches have tangent line with slope $0$ and to preserve it we must consider $\gamma =0$. So, the corresponding element in $\mathcal{A}$ is a composition of elements of $\mathcal{H}$ and $\tilde{\mathcal{A}}_1$.

If $\varphi \in \mathcal{P}$ has just two blocks with $\theta_1=0$ and $\theta_2=\infty$, then we must consider $T(\theta)=\frac{\delta\theta}{\alpha}$ to preserve the
slopes of the tangent lines and consequently the corresponding elements in $\mathcal{A}$ are obtained as a composition of elements of $\mathcal{H}$ and $\mathcal{A}_1$.

On the other hand, if $s\geq 3$ the branches of the first three blocks have tangent line with slope $0$, $\infty$ and $1$ respectively, and the only M\"{o}bius transformation that preserve them is the identity ($\gamma=\beta =0$ and $\alpha = \delta$). Consequently, the associated diffeomorphisms belong to $\mathcal{H}' \circ \mathcal{A}_1$.
\cqd

The standard strategy to solve the analytic equivalence problem is to find a  representative element in each orbit (a normal form), for instance, with a few number of parameters, in such a way that it is manageable to decide whether two normal forms are equivalent or not.

Notice that the action of any subgroup of $\mathcal{S}\times\mathcal{H}$ on an element $\varphi\in\mathcal{P}$ does not introduce or eliminate terms in the multigerm components. So, to obtain equivalent multigerms in the block form with short parameterizations it is natural to consider the $\mathcal{A}_1$-action or $\tilde{\mathcal{A}}_1$-action according to the number $s$ of blocks described in the previous lemma.

In a more general situation, concerning a Lie group action $\mathfrak{G}$ on an affine space $\mathbb{A}$, the {\it Complete Transversal Theorem} (see \cite{BruceKirkPlessis}) provides us a method for obtaining representative elements in the orbit $\mathfrak{G}(v)$ of $v\in \mathbb{A}$, if we have the description of the tangent space
$T\mathfrak{G}(v)$ to orbit $\mathfrak{G}(v)$ at $v$. More precisely, the following version will be
useful for our purposes:

\begin{theorem}{\sc (Complete Transversal Theorem (CTT))}\label{CTT} Let $\mathfrak{G}$ be a Lie group acting on an affine space $\mathbb{A}$ with underlying vector space $V$ and let $W$ be a subspace of $V$.
Suppose that $T\mathfrak{G}(v+w)=T\mathfrak{G}(v)$ for all $v \in \mathbb{A}$ and $w\in W$. If $v \in \mathbb{A}$ and  $W \subseteq T\mathfrak{G}(v)$, then $\mathfrak{G}(v+w)=\mathfrak{G}(v)$, for all $w\in W$, that is, the
$\mathfrak{G}$-orbits of $v+w$ and $v$ coincide.
\end{theorem}

Let us apply the above theorem in our context. To simplify the notation let $\mathcal{G}$ be one of the groups $\mathcal{A}_1$ or $\tilde{\mathcal{A}}_1$ and $I=\{1,\ldots ,r\}$.

The description of $T\mathcal{A}_1(\varphi)$ for $\varphi\in\mathcal{P}$ is classical (see for instance \cite{Wall}) and we do not need any effort to extend it to the group $\tilde{\mathcal{A}}_1$. Explicitly, we get
\begin{equation}
\label{TangSpace}
\begin{array}{c}
T{\mathcal{G}}(\varphi)=\left\{\left [\begin{array}{ccc}
x_1' \cdot \epsilon_1+ \varphi^*_1(\eta_1) & \ldots & x_r' \cdot\epsilon_r+ \varphi^*_r(\eta_1)\\
y_1' \cdot \epsilon_1+ \varphi^*_1(\eta_2) & \ldots & y_r' \cdot\epsilon_r+ \varphi^*_r(\eta_2)
\end{array}\right ]; \ \epsilon_i \in \langle t^2_i \rangle, i \in I \right\},\vspace{0.2cm}\\
\mbox{where} \ \left\{\begin{array}{ll}
\eta_1,\eta_2 \in \mathcal{M}^2,  & \mbox{if}\ \ \mathcal{G}=\mathcal{A}_1,\\
\eta_1 \in \langle X^2,Y \rangle,\ \eta_2 \in \mathcal{M}^2,
& \mbox{if}\ \ \mathcal{G}=\tilde{\mathcal{A}}_1.
\end{array}\right.
\end{array}
\end{equation}

Denoting by $B^k$ the set of $k$-jets of elements of a set $B$, we know that $\mathcal{G}^k$ is a unipotent Lie group that acts on $\mathcal{P}^k$ in the natural way: $$j^{k}(\rho_1,\ldots ,\rho_r,\sigma)\cdot j^{k}\varphi:=[j^{k}(\sigma\circ\varphi_1\circ\rho^{-1}_1),\ldots ,j^{k}(\sigma\circ\varphi_r\circ\rho^{-1}_r)].$$
The corresponding tangent space $T{\mathcal{G}}^k(j^k\varphi)$ is the set of $k$-jets of elements in $T{\mathcal{G}}(\varphi)$.

If $H^k$ denotes the $\mathbb{C}$-vector space
\begin{equation}
\label{SetHk}
\left\{ \left[\begin{array}{ccc}
c_1t_1^k  &  \dots & c_rt_r^k \\
d_1t_1^k  &  \dots &  d_rt_r^k \\
\end{array}\right] ; \ c_i,\ d_i\in\mathbb{C}\ \mbox{with}\ c_i=d_i=0\ \mbox{if}\ k \leq n_i,\ i\in I  \right\}
\end{equation}
then, we can prove, similarly as in Proposition 2 of \cite{HefezHernRod}, that $T{\mathcal{G}}^k(j^k\varphi+\zeta)=T{\mathcal{G}}^k(j^k\varphi)$ for any $j^k\varphi\in\mathcal{P}^k$ and $\zeta\in H^k$.

In this way, $\mathfrak{G}=\mathcal{G}^k$, $\mathbb{A}=\mathcal{P}^k$ and $W\subseteq H^k\cap T{\mathcal{G}}^k(j^k\varphi)$ fulfill the hypothesis of Theorem \ref{CTT}, that is, for any $w\in W$ we have that $j^k\varphi+w$ is $\mathcal{G}^k$-equivalent to $j^k\varphi$ and, consequently, there exists $\psi\in\mathcal{P}$ which is $\mathcal{A}$-equivalent to $\varphi$ with $j^k\psi=j^k\varphi+w$.

In the next two sections we will recognize elements in $H^k\cap T{\mathcal{G}}^k(j^k\varphi)$ using $\mathcal{G}$-invariants.

\section{Puiseux block form}
\label{section-puiseux}

In this section, we present a convenient parameterization for the components of the multigerm $\varphi=[\varphi_1,\ldots ,\varphi_r]\in\mathcal{P}$ preserving the block form, as in Definition \ref{blockform}. In addition, we exhibit some elimination criteria for terms of $\varphi$.

For simplicity, we denote $E_{ji}$ for $j=1,2$, $i\in I$, an element of the canonical basis of the complex matrices of order $2\times r$ and $\mathcal{G}$ denotes one of the groups $\mathcal{A}_1$ or $\tilde{\mathcal{A}}_1$.

Given a multigerm $\varphi=[\varphi_1,\ldots ,\varphi_r]\in\mathcal{P}$ with $\varphi_i(t_i)=(x_i(t_i),y_i(t_i))$, an element
$$\left[\begin{array}{ccc}
u_1 & \cdots & u_r \\
v_1 & \cdots & v_r
\end{array}\right]:=\sum_{i=1}^r (u_iE_{1i}+v_iE_{2i})$$ with $u_i,v_i\in\mathbb{C}\{t_i\}$ belongs to $T\mathcal{G}(\varphi)$ if and only if there exist
$\epsilon_i$ for $i=1,\ldots, r$ and $\eta=(\eta_1, \eta_2)$ satisfying (\ref{TangSpace}) such that
\begin{equation}\label{systemTang}\left\{\begin{array}{c}
u_i=x_i' \cdot \epsilon_i + \varphi^*_i(\eta_1) \\
\ v_i= y_i' \cdot \epsilon_i + \varphi^*_i(\eta_2). \end{array}\right.
\end{equation}

As we remarked in Section \ref{section-groupactions} describing $H^k\cap T\mathcal{G}^k(j^k\varphi)$ we can apply Theorem \ref{CTT} in order to obtain representative elements in a same $\mathcal{G}^k$-orbit.

Note that our analysis, up to this point, has taken into account a few $\mathcal{A}$-invariants, just the multiplicities of the branches and the behavior of the tangent cone. It is time to consider finer invariants as defined in Section \ref{section:AnalInvariants}.

Let $\varphi =[\varphi_1,\ldots ,\varphi_r]\in\mathcal{P}$ be a multigerm of plane curve in block form,  with local ring $\mathcal{O}$. The (fractional) ideal $\varphi^*(\mathcal{M}^2)$ of $\mathcal{O}$ plays a relevant role and provides us an elimination criterion for terms in $\varphi$.

\begin{lemma}\label{elim-gros} With the standard notations, if $k=\gamma+\sum_{{j\in I\atop j\neq i}}\nu_i(f_j)\in\nu_i(\mathcal{M}^2)$ with $\gamma\in\Gamma_i$, then there exists $\psi=[\psi_1,\ldots ,\psi_r]$  $\mathcal{G}$-equivalent to $\varphi$ such that $j^{k}\varphi_j=j^{k}\psi_j$ for all $j\in I\setminus\{i\}$ and $j^k\psi_i=j^{k-1}\varphi_i$.
\end{lemma}
\Dem By hypothesis, there exists an element $h\in\mathbb{C}\{X,Y\}$ with $h\cdot\prod_{{j\in I\atop j\neq i}}f_j\in \mathcal{M}^2$ and  $j^k(\varphi_i^*(h\cdot\prod_{{j\in I\atop j\neq i}}f_j))=t_i^k$.

Taking
$\epsilon_j=0$ for $j\in I$ and $\eta_l= \alpha_lh\cdot\prod_{{j\in I\atop j\neq i}}f_j$ for $l=1,2$ in (\ref{systemTang}), we get
$$\alpha_1 t_i^kE_{1i}+\alpha_2 t_i^kE_{2i}=\left[\begin{array}{ccccc}
0 & \cdots & \alpha_1t_i^{k} & \cdots & 0 \\
0 & \cdots & \alpha_2t_i^{k} & \cdots & 0
\end{array}\right] \in H^k\cap T\mathcal{G}^k(j^k \varphi).$$
By Theorem \ref{CTT}, we  obtain $\psi\in\mathcal{G}(\varphi)$ with $j^k\psi=j^k\varphi+\alpha_1 t_i^kE_{1i}+\alpha_2 t_i^kE_{2i}$ and for a convenient choice of $\alpha_1,\alpha_2\in\mathbb{C}$ we have $j^{k}\psi_j=j^{k}\varphi_j$ for all $j\in I\setminus\{i\}$ and $j^k\psi_i=j^{k-1}\varphi_i$.
\cqd

\begin{remark}\label{mudanca1} Notice that the description of tangent vector $\alpha_1 t_i^kE_{1i}+\alpha_2 t_i^kE_{2i}\in H^k\cap T\mathcal{G}^k(j^k \varphi),$ $k=\gamma+\sum_{{j\in I\atop j\neq i}}\nu_i(f_j)\in\nu_i(\mathcal{M}^2)$ with $\gamma\in\Gamma_i$ in the previous lemma gives us a clue to obtain elements $(\rho_1,\ldots ,\rho_r,\sigma)\in\mathcal{G}$ such that $(\rho_1,\ldots ,\rho_r,\sigma)\cdot\varphi=\psi$ with $j^{k}\psi_j=j^{k}\varphi_j$ for all $j\in I\setminus\{i\}$ and $j^k\psi_i=j^{k-1}\varphi_i$. In fact, with the above notations, it is sufficient to consider $\rho_j(t_j)=t_j$ for $j\in I$ and $\sigma (X,Y)=(X-\alpha_1\cdot h\cdot\prod_{{j\in I\atop j\neq i}}f_j,\ Y-\alpha_2\cdot h\cdot\prod_{{j\in I\atop j\neq i}}f_j)$ with appropriate $\alpha_1,\alpha_2\in\mathbb{C}$.
\end{remark}

As a consequence of the above lemma we obtain an estimate for the finite determinacy of a given $\varphi\in\mathcal{P}$.

Denoting $(\mathcal{O}:\overline{\mathcal{O}})^c$ the contraction of $(\mathcal{O}:\overline{\mathcal{O}})$ by $\varphi^*$ we have the following result.

\begin{proposition}\label{trunc} If $d_i$ is the conductor of the set $\nu_i\left (\mathcal{M}^2\cap (\mathcal{O}:\overline{\mathcal{O}})^c\right )$ for $i=1, \ldots, r$, then
    $\varphi\underset{\mathcal{G}}{\sim}[j^{d_1-1}\varphi_1,\ldots ,j^{d_r-1}\varphi_r].$
\end{proposition}
\Dem For each $i\in I$ we put $\varphi_i-j^{d_i-1}\varphi_i=(t_i^{d_i}u_{i1},t_i^{d_i}u_{i2})$ with $u_{il}\in\mathbb{C}\{t_i\}$ for $l=1,2$. As $d_i$ is the conductor of $\nu_i\left (\mathcal{M}^2\cap (\mathcal{O}:\overline{\mathcal{O}})^c\right )$ there exists $h_{il}\cdot\prod_{j\in I\atop j\neq i}f_j\in \mathcal{M}^2\cap (\mathcal{O}:\overline{\mathcal{O}})^c$ such that $\varphi_i^*(h_{il}\cdot\prod_{j\in I\atop j\neq i}f_j)=t_i^{d_i}u_{il}$ for $l=1,2$.

Now, taking $(\rho_1,\ldots, \rho_r,\sigma)\in \mathcal{G}$, where $\rho_j(t_j)=t_j$ for $j\in I$ and $$\sigma (X,Y)=\left (X-\sum_{i\in I} h_{i1}\cdot\prod_{j\in I\atop j\neq i}f_j,\ Y-\sum_{i\in I} h_{i2}\cdot\prod_{{j\in I\atop j\neq i}}f_j\right ),$$ we obtain  $\psi=(\rho_1,\ldots ,\rho_r,\sigma)\cdot\varphi=[j^{d_1-1}\varphi_1,\ldots ,j^{d_r-1}\varphi_r]$.
\cqd

Since $\kappa=(\kappa_1, \ldots, \kappa_r)$ is the conductor of the semiring $\Gamma$, the above mentioned integer $d_i$ satisfies $d_i\geq \kappa_i$ with equality if
$(\mathcal{O}:\overline{\mathcal{O}})^c\subseteq \mathcal{M}^2$.
More precisely, we have the following:

\begin{proposition}\label{di} Let $\varphi=[\varphi_1, \ldots, \varphi_r]\in\mathcal{P}$ be a block form multigerm and $s$ the number of blocks. For each $i \in I$, the integer $d_i$ in the last proposition   satisfies
    $$d_i=\left \{ \begin{array}{ll}
    \kappa_i+2 & \mbox{if}\ \ r=1\ \mbox{and}\ n_i\leq 2; \\
    \kappa_i+1 & \mbox{if}\ \ r=s=2\ \mbox{and}\ n_1=n_2=1; \\
    \kappa_i & \mbox{otherwise}.
    \end{array}\right.$$
\end{proposition}
\Dem Let us consider the cases $r=1,\ r=2$ and $r\geq 3$ separately:

{\bf Case $r=1$:}

Let $\varphi_1$ be a parameterization of a branch with multiplicity
$n_1$. As we remarked in Section \ref{section:AnalInvariants}, $\kappa_1$ is the Milnor number $\mu_1$.

If $n_1=1$, then $\mu_1=0$ and $\mathcal{O}=\overline{\mathcal{O}}$. So, $d_1=\min\ \nu_i\left (\mathcal{M}^2\cap (\mathcal{O}:\overline{\mathcal{O}})^c\right )=2=\mu_1+2=\kappa_1+2$.

For $n_1>1$ we have $\overline{\mathcal{O}}\neq\mathcal{O}$ and
$(\mathcal{O}:\overline{\mathcal{O}})^c \subseteq \mathcal{M}$. Setting
$m=\min(\Gamma\setminus n_1\mathbb{N})$, we have $$\nu_1
(\mathcal{M}\setminus\mathcal{M}^2)=\left \{n_1, 2n_1,\ldots , \left [\frac{m}{n_1}\right ]n_1, m\right \}.$$

If $n_1=2$, then $\Gamma =\langle 2,m\rangle=2\mathbb{N}+m\mathbb{N}$ and $\mu_1=m-1$, which implies that $\nu_1((\mathcal{O}:\overline{\mathcal{O}})^c)=\{\gamma\in\mathbb{N};\ \gamma\geq m-1\}$. So, the conductor $d_1$ of $\nu_1\left (\mathcal{M}^2\cap (\mathcal{O}:\overline{\mathcal{O}})^c\right )$ satisfies $d_1=m+1=\mu_1+2$. On the other hand, for $n_1>2$ we have $\mu_1>m$, $(\mathcal{O}:\overline{\mathcal{O}})^c\subseteq \mathcal{M}^2$ and consequently $d_1=\kappa_1$.

{\bf Case $r=2$:}

Consider a plane curve $\mathcal{C}$ given by $f=f_1\cdot f_2$
with corresponding multigerm $[\varphi_1,\varphi_2]$.

As $(\mathcal{O}:\overline{\mathcal{O}})=(\mathcal{O}_1:\overline{\mathcal{O}_1})\cdot \varphi_1^*(f_2)\oplus (\mathcal{O}_2:\overline{\mathcal{O}_2})\cdot \varphi_2^*(f_1)$ and $f_1,f_2\in \mathcal{M}$, it follows that if $n_i>1$ for some $i\in\{1,2\}$, then $(\mathcal{O}:\overline{\mathcal{O}})^c\subseteq \mathcal{M}^2$ and $d_i=\kappa_i$.

Let us consider $n_1=n_2=1$. In this case we have $\kappa=(\kappa_1,\kappa_2)=(\nu_1(f_2),\nu_2(f_1))$. If $s=1$, then $\varphi_1,\varphi_2\in B_1$ and $\nu_1(f_2)=\nu_2(f_1)>1$. Hence,
$(\mathcal{O}:\overline{\mathcal{O}})^c\subseteq \mathcal{M}^2$ and $d_i=\kappa_i$ for $i=1,2$. If $s=2$, then $\varphi_1\in B_1$ and $\varphi_2\in B_2$. So, $\nu_1(f_2)=\nu_2(f_1)=1$ and $d_i=\min\ \nu_i\left (\mathcal{M}^2\cap (\mathcal{O}:\overline{\mathcal{O}})^c\right )=2=\kappa_i+1$ for $i=1,2$.

{\bf Case $r\geq 3$:}

Now, we consider a plane curve $\mathcal{C}$ given by $f=f_1\cdot
\ldots\cdot f_r$ with $f_i\in \mathcal{M}$.

As $(\mathcal{O}:\overline{\mathcal{O}})=\bigoplus_{i\in
    I}(\mathcal{O}_i:\overline{\mathcal{O}_i})\cdot \prod_{j\in I\atop
    j\neq i}\varphi_i^*(f_j)$, we get
$(\mathcal{O}:\overline{\mathcal{O}})^c\subseteq \mathcal{M}^2$ and
$d_i=\kappa_i$ for $i\in I$. \cqd

Another direct application of Theorem \ref{CTT} is obtained as consequence of the following lemma:

\begin{lemma}
    \label{LemmaElemTang}
    Let $\varphi=[\varphi_1, \ldots, \varphi_r]\in\mathcal{P}$ be a block form multigerm associated to a curve $\mathcal{C}$ and
    $\theta_i=\frac{b_{i}}{a_{i}}\in \overline{\mathbb{C}}$ the slope of the tangent line to $\mathcal{C}_i$ at the origin. For each $i\in I$, if $k> n_i$, then
    $$t_i^k(a_{i}E_{1i}+b_{i}E_{2i})\in H^k\cap T\mathcal{G}^k(j^k\varphi).$$
\end{lemma}
\Dem As $\varphi_i(t_i)=(a_{i}t_i^{n_i}+\ldots,
b_{i}t_i^{n_i}+\ldots)$ and  $k>n_i$, taking
$\epsilon_i=\frac{1}{n_i}t_i^{k-n_i+1}\in \langle t_i^2
\rangle$, $\epsilon_l=0$ for $l\in I\setminus\{i\}$ and
$\eta_1=\eta_2=0$ in (\ref{systemTang}), we obtain
$$t_i^k(a_{i}E_{1i}+b_{i}E_{2i})=\left[\begin{array}{ccccc}
0 & \cdots & a_{i}t_i^{k}& \cdots & 0 \\
0 & \cdots & b_{i}t_i^{k}& \cdots & 0
\end{array}\right] \in H^k\cap T\mathcal{G}^k(j^k \varphi).$$\cqd

\vspace{-0.5cm}
By the above lemma and Theorem \ref{CTT}, for every $k>n_i$ there exists a block form multigerm $\psi\in\mathcal{G} (\varphi)$ with $j^k\psi=j^k\varphi+\alpha t_i^k(a_{i}E_{1i}+b_{i}E_{2i})$ for any $\alpha\in\mathbb{C}$. So, if $\varphi_i=(x_i,y_i)\in B_1$, that is, $a_i=1$ and $b_i=0$, then
we can choose $\alpha$ in such way that $j^k\psi_i=(j^{k-1}x_i,j^{k}y_i)$, that is, we eliminate the $k$-order term of the first component of $\varphi_i$. On the other hand for $\varphi_i\in B_2$ there exists $\alpha$ such that $j^k\psi_i=(j^kx_i,j^{k-1}y_i)$. Finally, for $\varphi_i\in B_j$ with $j>2$, as $a_i\cdot b_i\neq 0$, we can choose $\alpha$ to obtain  $j^k\psi_i=(j^{k-1}x_i,j^{k}y_i)$ or $j^k\psi_i=(j^kx_i,j^{k-1}y_i)$.

Similarly to Remark \ref{mudanca1}, as consequence of the previous lemma, we can exhibit an element of $\mathcal{G}$ to perform the proposed action. In fact, with the above notations, it is sufficient to consider $\rho_i^{-1}(t_i)=t_i+\alpha\epsilon_i$ with $\alpha\in\mathbb{C}$, $\epsilon_i\in\langle t_i^2\rangle$ given as in the proof of Lemma \ref{LemmaElemTang}, $\epsilon_j=0$ for $j\in I\setminus\{i\}$ and $\sigma (X,Y)=(X,Y)$.

Thus we recover the classical Puiseux expansion for plane curves.

\begin{proposition}\label{prop-puiseux-block}
    Any multigerm $\varphi=[\varphi_1, \ldots, \varphi_r]\in\mathcal{P}$ is $\mathcal{G}$-equivalent to a block form multigerm $\psi=[\psi_1,\ldots ,\psi_r]$ with
    \begin{equation}\label{ParamBlocks}
    \psi_i(t_i)=\left\{\begin{array}{ll}
    \left (t_i^{n_i},\ \theta_lt_i^{n_i}+\displaystyle\sum_{j>n_i}a_{ij}t_i^j\right ) & \mbox{for}\ \psi_i\in B_l,\ l\neq 2;\vspace{0.1cm}\\
    \left (\displaystyle\sum_{j>n_i}a_{ij}t_i^j,\ t_i^{n_i}\right ) & \mbox{for}\ \psi_i\in B_2;\vspace{0.1cm}
    \end{array}\right.
    \end{equation}
    with $\theta_1=0$, $\theta_3=1$ and $\theta_k \neq \theta_l$ for $1\leq k, l\leq r$ and $k\neq l$ .
\end{proposition}
\Dem We may suppose that $\varphi$ is given by a block form multigerm. For each $i\in I$ and for every $k>n_i$ we consider elements in $\mathcal{G}$ to apply the process described after Lemma \ref{LemmaElemTang}. If $\varphi_i=(x_i,y_i)\in B_l$ with $l\geq 3$, then we choose $j^k\psi_i=(j^{k-1}x_i,j^{k}y_i)$.

Repeating this process for all $k$ satisfying $n_i<k<d_i$ we obtain
$j^{d_i}\psi_i=(t_i^{n_i},j^{d_i}y_i)$ for $\psi_i\not\in B_2$ and
 $j^{d_i}\psi_i=(j^{d_i}x_i,t_i^{n_i})$ for $\psi_i\in B_2$. The result follows from Proposition \ref{trunc}.
\cqd

We call the multigerm as highlighted in the above result the \textbf{Puiseux block form} of $\varphi$.

From above results, we obtain the well
known representative in the analytical class for the particular
cases of Proposition \ref{di}. More explicitly, if $\varphi$ is a Puiseux
block form, then:
$$\varphi\underset{\mathcal{A}}{\sim} (t,0)\ \mbox{if}\ r=n_1=1;$$
$$\varphi\underset{\mathcal{A}}{\sim} (t^2,t^{m})\ \mbox{if}\ r=1,\ n_1=2\ \mbox{and}\ m=\min(\Gamma\setminus 2\mathbb{N});$$
$$\varphi\underset{\mathcal{A}}{\sim} [\psi_1,\psi_2];\ \psi_1=(t_1,0), \psi_2=(0,t_2)\ \mbox{if}\ r=s=2,\ n_1=n_2=1.$$

Notice that if $r=1$, $n_1=2$ and $m=\min(\Gamma\setminus 2\mathbb{N})$ the previous results provide us that $\varphi\underset{\mathcal{A}}{\sim} (t^2,at^{m})$ with $a\neq 0$. But, considering $\rho(t)=t$ and $\sigma (X,Y)=(X,a^{-1}Y)$ we obtain the above equivalence.

In what follows, we will consider the non exceptional above cases and consequently we can suppose that the finite determinacy order for each $\varphi_i$ of a Puiseux block form multigerm is $\kappa_i$.

\section{$\mathcal{G}$-Normal forms and Analytic Equivalence}
\label{section-normalform}

In Section \ref{section-puiseux}, we apply the Complete Transversal Theorem (Theorem \ref{CTT}) to reduce a multigerm $\varphi\in\mathcal{P}$ to its Puiseux block form taking into account only a few topological invariants. Although such forms are a shorter way of presenting a multigerm it is not easy to decide if two Puiseux block form multigerms correspond or not to analytically equivalent plane curves.

Our goal in this section is to partition a topological class into a finite number of strata so that in each stratum we have a constant $\mathcal{A}$-invariant and every element belonging to it admits a particular Puiseux block form, that we call a {\it normal form}, which will allow us to distinguish them from the analytical equivalence viewpoint.

As before, according to Lemma \ref{subgroups}, $\mathcal{G}$ denotes $\mathcal{A}_1$ if the Puiseux block form of $\varphi\in\mathcal{P}$ has at least two blocks and $\tilde{\mathcal{A}}_1$ if there is a single block.

Recall that as a consequence of Lemma \ref{LemmaElemTang},  $\{\sum_{\varphi_i\not\in B_2}\mathbb{C}\cdot t_i^kE_{1i}+\sum_{\varphi_i\in B_2}\mathbb{C}\cdot t_i^kE_{2i}\}\cap H^k \subseteq T\mathcal{G}^k(j^k \varphi)$ for every $k$. Thus, in order to obtain a Puiseux block form multigerm with a smaller number of terms, it is sufficient to describe elements in the set $D^k:=\{\sum_{\varphi_i\not\in B_2}\mathbb{C}\cdot t_i^kE_{2i}+\sum_{\varphi_i\in B_2}\mathbb{C}\cdot t_i^kE_{1i}\}\cap H^k$ that belong to $T\mathcal{G}^k(j^k \varphi)$. In fact, if $\zeta\in D^k\cap T\mathcal{G}^k(j^k\varphi)$ then $\varphi$ is $\mathcal{G}$-equivalent to a Puiseux block form  multigerm such that its $k$-jet is $j^k\varphi+\zeta$, that is, there exists a multigerm $\psi$, $\mathcal{A}$-equivalent to $\varphi$, with $j^{k-1}\psi=j^{k-1}\varphi$ and, by adjusting coefficients in $\zeta$, the element $j^k\psi$ has a smaller number of nonzero terms than $j^k\varphi$.

The key to obtain our normal form is the connection between
 $D^k\cap T\mathcal{G}^k(j^k\varphi)$
and the set of values of a particular fractional ideal as we describe in the sequel.

Considering $\varphi \in \mathcal{P}$ and
$$
\Omega_{\mathcal{A}_1}=\left \{\eta_2dX-\eta_1dY; \ \eta_1,\eta_2 \in \mathcal{M}^2\right \}, \ \ \Omega_{\tilde{\mathcal{A}}_1}=\left \{\eta_2dX-\eta_1dY; \ \eta_1\in\langle X^2,Y\rangle\ \mbox{and}\ \eta_2\in \mathcal{M}^2\right \},$$
we get the fractional ideal $$\mathcal{I}_{\mathcal{G}}=\left \{\left (\frac{\varphi^{*}_{1}(\omega)}{n_1t_1^{n_1}},\ldots , \frac{\varphi^{*}_{r}(\omega)}{n_rt_r^{n_r}}\right );\ \omega\in\Omega_{\mathcal{G}}\right \}$$ of $\mathcal{O}$ according to the group $\mathcal{G}$.

We define $\Lambda_{\mathcal{G}}:=\{\nu({\rm w});\ {\rm w} \in \mathcal{I}_{\mathcal{G}}\}\subset \Lambda-n \subseteq \overline{\mathbb{N}}^r$, where $\Lambda$ is the $\mathcal{A}$-invariant defined in (\ref{LambdaInvariant})  and $n:=(n_1,\ldots, n_r)$ with $n_i$ the multiplicity of the branch $\varphi_i$.

Notice that the conductor $\varrho'=(\varrho_1',\ldots, \varrho_r')$ of $\Lambda_{\mathcal{G}}$ satisfies $\varrho_i' \leq \varrho_i-n_i < \kappa_i$, where $(\varrho_1, \ldots, \varrho_r)$ and $(\kappa_1, \ldots, \kappa_r)$ are the conductors of $\Lambda$ and $\Gamma$, respectively. As a consequence of Remark \ref{Remark-Maximal-Box}, for each fixed value semiring $\Gamma$ there are only finitely many $\mathcal{A}$-invariants $\Lambda_{\mathcal{G}}$ that can be computed by the algorithm presented in Theorem 16, \cite{carvalho2}.

In what follows we will consider the set $\Sigma_{\Gamma,\Lambda_{\mathcal{G}}}$ of all Puiseux block form multigerms in $\mathcal{P}$, as given in (\ref{ParamBlocks}), with semiring $\Gamma$ and $\Lambda_{\mathcal{G}}$ fixed.

\begin{proposition}
\label{DifKTg} Given $\varphi=[\varphi_1, \ldots, \varphi_r]\in \Sigma_{\Gamma,\Lambda_{\mathcal{G}}}$ and $k\in\mathbb{N}^*$ there exists ${\rm w}=({\rm w}_1,\ldots ,{\rm w}_r) \in\mathcal{I}_{\mathcal{G}}$ with $j^k{\rm w}_i=\alpha_it_i^{k}$, $\alpha\in\mathbb{C}$ if and only if
 $\sum_{\varphi_i\not\in B_2} \alpha_it_i^{k}E_{2i}-\sum_{\varphi_i\in B_{2}} \alpha_it_i^{k}E_{1i} \in T{\mathcal{G}}^k(j^k\varphi)$.
\end{proposition}
\Dem Let ${\rm w}=({\rm w}_1,\ldots ,{\rm w}_r)\in\mathcal{I}_{\mathcal{G}}$ with $j^k{\rm w}_i=j^k\left (\frac{\varphi^{*}_{i}(\omega)}{n_it_i^{n_i}}\right )=\alpha_it_i^{k}$ for $i\in I$ and some $\omega=\eta_2dX-\eta_1dY\in\Omega_{\mathcal{G}}$.

As
$\epsilon_i = -\frac{\varphi^*_i(\eta_1)}{n_it_i^{n_i-1}}\in\langle t_i^2\rangle$ for $\varphi_i \not\in B_2$, $\epsilon_i =- \frac{\varphi^*_i(\eta_2)}{n_it_i^{n_i-1}}\in\langle t_i^2\rangle$ for
$\varphi_i \in B_{2}$,  $\eta_1$ and $\eta_2$ satisfying (\ref{systemTang}) we have the desired element in $T{\mathcal{G}}^k(j^k\varphi)$.

Conversely, if  $\sum_{\varphi_i\not\in B_2} \alpha_it_i^{k}E_{2i}-\sum_{\varphi_i\in B_{2}} \alpha_it_i^{k}E_{1i} \in T{\mathcal{G}}^k(j^k\varphi)$, then there exist $\epsilon_i \in \langle t_i^2 \rangle,$ $i\in I$, $\eta_1$ and $\eta_2$ as given in (\ref{TangSpace}) such that
\begin{equation}\label{alpha}
\left\{\begin{array}{l}
x_i' \cdot \epsilon_i + \varphi^*_i(\eta_1)=P_i \\
y_i' \cdot \epsilon_i + \varphi^*_i(\eta_2)=\alpha_it_i^k+Q_i\end{array}\right. \mbox{if}\ \varphi_i\not\in B_2,\ \ \  \left\{\begin{array}{l}
x_i' \cdot \epsilon_i + \varphi^*_i(\eta_1)=-\alpha_it_i^k+P_i \\
y_i' \cdot \epsilon_i + \varphi^*_i(\eta_2)=Q_i\end{array}\right. \mbox{if}\ \varphi_i\in B_2
\end{equation}
with $\nu_i(P_i)>k$ and $\nu_i(Q_i)>k$ for $i\in I$. Considering
$\omega = \eta_2 dX - \eta_1 dY \in \Omega_{\mathcal{G}}$ we provide the element ${\rm w}\in \mathcal{I}_{\mathcal{G}}$ satisfying the proposition. \cqd

Notice that if $k\leq n_i$ then, by definition of $\mathcal{I}_{\mathcal{G}}$, $\alpha_i$ should be zero in (\ref{alpha}).

In order to simplify the notation, from now on given a non empty subset $J\subseteq I$, we denote the $J$-fiber of ${\underline k}=(k, \dots, k)\in\mathbb{N}^r$ with respect to $\Lambda_{\mathcal{G}}$ by $F_J(\underline{k})$. The following corollary is immediate.

\begin{corollary}
\label{CorDifTg}
We have $F_J(\underline{k})\neq\emptyset$ if and only if there exist $\alpha_i\neq 0$ for every $i\in J$ and $\alpha_i= 0$ for every $i\in I\setminus J$ such that $\sum_{\varphi_i\not\in B_2} \alpha_it_i^kE_{2i}-\sum_{\varphi_i\in B_{2}} \alpha_it_i^kE_{1i} \in T{\mathcal{G}}^k(j^k\varphi)$.
\end{corollary}

The next proposition describes subspaces of $D^k\cap T{\mathcal{G}}^k(j^k\varphi)$ in order to choose terms which can be eliminated by $\mathcal{G}$-action. Such subspaces are related to the fibers of $\underline{k}$ with respect to $\Lambda_{\mathcal{G}}$.

\begin{proposition}
\label{PropElim} Given $k\in \mathbb{N}^*$ there exists a
$d$-dimensional subspace in $D^k\cap T{\mathcal{G}}^k(j^k\varphi)$ if and only
if there exists $L=\{l_1, \ldots, l_{d}\}\subseteq I$ satisfying the
following condition:
\begin{center}
For each $l\in L$, there exists $J_l \subseteq I$ with $l \in J_l \setminus \bigcup_{i\in L\setminus\{l\}}J_i$, and $F_{J_l}(\underline{k})\neq \emptyset\ \ \ \ \ (\star).$
\end{center}
\end{proposition}
\Dem We will show that there exist $v_1,\ldots ,v_d\in D^k\cap T{\mathcal{G}}^k(j^k\varphi)$ linearly independent if and only if there exists $L=\{l_1, \ldots, l_{d}\}\subseteq I$ satisfying $(\star)$.

Let $\{v_i:=\sum_{\varphi_j\not\in B_2} \alpha_{ij}t_j^kE_{2j}+\sum_{\varphi_j\in B_{2}} \alpha_{ij}t_j^kE_{1j}\in D^k\cap T{\mathcal{G}}^k(j^k\varphi);\ 1\leq i\leq d\}$ be a set of linearly
independent vectors. Then there exists $L=\{l_1, \ldots, l_{d}\}\subseteq I$
with $[\alpha_{il_j}]_{1\leq i,j\leq d}\in Gl_d(\mathbb{C})$.

In this way, for each $1\leq i\leq d$ there exists $w_i=\sum_{\varphi_j\not\in B_2}
\beta_{ij}t_j^kE_{2j}+\sum_{\varphi_j\in B_{2}} \beta_{ij}t_j^kE_{1j}\in
span\{v_1,\ldots ,v_d\}\subseteq D^k\cap T{\mathcal{G}}^k(j^k\varphi)$ with
$\beta_{il_i}=1$ and $\beta_{ij}=0$ for $j\in L\setminus\{l_i\}$. Taking the set
$J_{l_i}=\{j\in I;\ \beta_{ij}\neq 0\}$, the previous corollary ensures the
condition ($\star$).

On the other hand, suppose that there is $L=\{l_1, \ldots,
l_{d}\}\subseteq I$ satisfying the condition $(\star)$. So, by Corollary
\ref{CorDifTg}, for each $i=1,\ldots ,d$, there exists
$v_i=\sum_{\varphi_j\not\in B_2}
\alpha_{ij}t_j^kE_{2j}+\sum_{\varphi_j\in B_{2}} \alpha_{ij}t_j^kE_{1j}\in
D^k\cap T{\mathcal{G}}^k(j^k\varphi)$ such that $\alpha_{ij}\neq 0$ if $j\in J_{l_i}$
and $\alpha_{ij}=0$ for $j\in I\setminus J_{l_i}$. As the matrix $[\alpha_{il_j}]_{1\leq i,j\leq d}$ is invertible, it
follows that $v_1,\ldots ,v_d$ are linearly independent. \cqd

\begin{remark}\label{lk-diff} Notice that, by previous result, if $L=\{l_1, \ldots,
l_{d}\}\subseteq I$ satisfies the condition $(\star)$ then there exist $v_1,\ldots ,v_d\in
D^k\cap T{\mathcal{G}}^k(j^k\varphi)$ linearly independent with
$v_i=\sum_{\varphi_j\not\in B_2}
\alpha_{ij}t_j^kE_{2j}+\sum_{\varphi_j\in B_{2}} \alpha_{ij}t_j^kE_{1j}$ such that $\alpha_{ij}\neq 0$ if $j\in J_{l_i}$ and $\alpha_{ij}=0$ for $j\in I\setminus J_{l_i}$. In particular, for any $b_{l_1},\ldots ,b_{l_d}\in\mathbb{C}$ there exists $v\in span\{v_1,\ldots ,v_d\}$ with $v=\sum_{\varphi_j\not\in B_2}
b_jt_j^kE_{2j}+\sum_{\varphi_j\in B_{2}} b_jt_j^kE_{1j}\in
D^k\cap T{\mathcal{G}}^k(j^k\varphi)$, or equivalently, by Proposition \ref{DifKTg}, there exists  ${\rm w}=({\rm w}_1,\ldots, {\rm w}_r)\in \mathcal{I}_{\mathcal{G}}$ with $j^k{\rm w}_j=b_jt_j^k$ for any $b_{l_1},\ldots ,b_{l_d}\in\mathbb{C}$.
\end{remark}

As an immediate consequence we can establish an elimination criterion based on
the set $\Lambda_{\mathcal{G}}$.

\begin{corollary}
\label{CorElim}
  Let $k\in \mathbb{N}^*$, $\varphi\in\mathcal{P}$ as in (\ref{ParamBlocks}) and $L=\{l_1, \ldots,
  l_{d}\}\subseteq I$ satisfying
  the condition $(\star)$. Then there exists $\psi \in \mathcal{P}$ such that $\psi
\underset{\mathcal{G}}{\sim}\varphi$, $j^{k-1}\psi=j^{k-1}\varphi$
   and $j^{k}\psi_{l}=j^{k-1}\varphi_{l}$ for all $l\in L$.
\end{corollary}
\Dem Since $\varphi$ is given as in (\ref{ParamBlocks}) we have that $j^k\varphi=j^{k-1}\varphi + w$, with $w=\sum_{\varphi_j\not\in B_2}
\beta_{j}t_j^kE_{2j}+\sum_{\varphi_j\in B_{2}} \beta_{j}t_j^kE_{1j}\in D^k$.
By the above remark, there exists $v=\sum_{\varphi_j\not\in B_2}
b_jt_j^kE_{2j}+\sum_{\varphi_j\in B_{2}} b_jt_j^kE_{1j}\in
D^k\cap T{\mathcal{G}}^k(j^k\varphi)$ with $b_l=-\beta_l$ for all $l\in L$. By Theorem \ref{CTT}, it follows that $\varphi$ is $\mathcal{G}$-equivalent to some $\psi\in\mathcal{P}$ with $j^{k}\psi=j^{k}\varphi +v=j^{k-1}\varphi+w+v$, that is, $j^{k-1}\psi=j^{k-1}\varphi$
and $j^{k}\psi_{l}=j^{k-1}\varphi_{l}$ for all $l\in L$.
\cqd

Given $k>0$ the parameter elimination method established in the above result depends upon a set $L$ satisfying the condition $(\star)$ in the Proposition \ref{PropElim}. The aim is to eliminate terms of order $k$ from the greater number of branches prioritizing the first ones. We can formalize this by defining the injective map $S:\{L;\ L\subseteq I\}\to \mathbb{N}^r$ given by $S(L)=(z_1, \ldots, z_r)$ where
$$z_i=\left\{\begin{array}{l}
1, \ \mbox{if}\ i \in L\\
0, \ \mbox{if}\ i \not\in L.\\
\end{array}\right.$$

\begin{definition}\label{Lk}
For each fixed $k$, we consider the set of indices $L_k\subseteq \{1,\ldots ,r\}$ such that $S(L_k)=max_{GLex}\{S(L);\ L\ \mbox{satisfying}\ (\star)\ \mbox{in Proposition \ref{PropElim}}\}$, where $max_{GLex}$ means the maximum with respect to graded lexicographic order.
\end{definition}

\begin{remark}\label{propriedade-max}
Let $\varrho$ be the conductor of $\Lambda_{\mathcal{G}}$ and $\underline{k}\in\varrho+\mathbb{N}^r$. We have that $F_{\{i\}}(\underline{k})\neq\emptyset$ for all $i\in I$, consequently $L_k=I$.

If $\underline{k}$ is an absolute maximal of $\Lambda_{\mathcal{G}}$, then $F_I(\underline{k})\neq\emptyset$ and $F_{J}(\underline{k})=\emptyset$ for all $J\subsetneq I$. In this case, $L_k=\{1\}$.

If $\underline{k}$ is a relative maximal of $\Lambda_{\mathcal{G}}$, then $F_{\{i\}}(\underline{k})=\emptyset$ for all $i\in I$ and $F_J(\underline{k})\neq\emptyset$ for all $J\subseteq I$ with $\sharp J>1$. Notice that for any $i_0\in I$ the set $L=I\setminus\{i_0\}$ and $J_l=\{l,i_0\}$ for all $l\in L$ satisfy the condition $(\star)$. Then, in this case we have $L_k=I\setminus\{r\}$.
\end{remark}

Now we can establish one of the main results of this paper:

\begin{theorem}($\mathcal{G}$-{\sc normal form})\label{A1normal-form}
With the above notation any $\varphi=[\varphi_1,\ldots, \varphi_r]\in \Sigma_{\Gamma,\Lambda_{\mathcal{G}}}$ is $\mathcal{G}$-equivalent to $\psi=[\psi_1, \ldots, \psi_r]\in\mathcal{P}$ such that
\begin{equation}\label{normalBlocks}
\psi_i(t_i)=\left\{\begin{array}{cl}
\left (t_i^{n_i},\ \theta_lt_i^{n_i}+\displaystyle\sum_{j>n_i\atop i\not\in L_j}a_{ij}t_i^j\right ) & \mbox{for}\ \psi_i\in B_l,\ l\neq 2;\vspace{0.1cm}\\
\left (\displaystyle\sum_{j>n_i\atop i\not\in L_j}a_{ij}t_i^j,\ t_i^{n_i}\right ) & \mbox{for}\ \psi_i\in B_2, \\
\end{array}\right.
\end{equation}
where $L_j$ is given in Definition \ref{Lk}, $\theta_1=0$, $\theta_3=1$ and $\theta_k \neq \theta_l$ for $1\leq k, l\leq r$ with $k \neq l$. Moreover, if $\psi'=[\psi'_1, \ldots, \psi'_r]\in\mathcal{P}$ is another multigerm $\mathcal{G}$-equivalent to $\varphi$ as (\ref{normalBlocks}) then $\psi'=\psi$, that is, the $\mathcal{G}$-normal form of $\varphi$ is unique.
\end{theorem}
\Dem By Proposition \ref{prop-puiseux-block} we can consider $\varphi\in\mathcal{P}$ in its Puiseux block form. For all $k>0$ taking $L_k$ as above defined and applying the Corollary \ref{CorElim} we obtain a $\mathcal{G}$-normal form $\psi$ as (\ref{normalBlocks}). Remark that by Proposition \ref{trunc} it is sufficient to consider $k<max\{d_i;\ i\in I\}$  where $d_i$ is given in Proposition \ref{di}.

Now suppose $\psi'=[\psi'_1, \ldots, \psi'_r]\in\mathcal{P}$ as (\ref{normalBlocks}), $\mathcal{G}$-equivalent to $\varphi$ and $\psi'\neq \psi$. Let us consider $k=\min\{l;\ j^{l}\psi\neq j^{l}\psi'\}$.

By Proposition \ref{PropElim} there exists $\{v_i\in D^k\cap T\mathcal{G}^k(j^k\psi);\ l_i\in L_k\}$ linearly independent (possibly empty) and subsets $J_{l_i}\subseteq I$ satisfying the condition $(\star)$ with
$v_i=\sum_{\psi_j\not\in B_2}
\alpha_{ij}t_j^kE_{2j}+\sum_{\psi_j\in B_{2}} \alpha_{ij}t_j^kE_{1j}$ such that $\alpha_{ij}\neq 0$ if $j\in J_{l_i}$ and $\alpha_{ij}=0$ for $j\in I\setminus J_{l_i}$.

Let $\mathcal{F}^k\subseteq \mathcal{G}^k$ be the subgroup that leave invariant the affine space $N^k=j^{k-1}\psi+V$ where $$V=\left \{\sum_{\psi_j\not\in B_2}
\beta_{j}t_j^kE_{2j}+\sum_{\psi_j\in B_{2}} \beta_{j}t_j^kE_{1j}\in D^k;\ \beta_j\in\mathbb{C}\ \mbox{and}\ \beta_l=0\ \mbox{for all}\ l\in L_k\right \}.$$

As the ground field is $\mathbb{C}$ and $\mathcal{F}^k$ is unipotent it follows that $\mathcal{F}^k$ is a connected group. So, the orbit $\mathcal{F}^k(j^k\psi)\subseteq N^k\subseteq \mathcal{G}^k(j^k\psi)$ is connected and it contains $j^k\psi$ and $j^k\psi'$. Let
$w:=\sum_{\psi_j\not\in B_2}
\gamma_{j}t_j^kE_{2j}+\sum_{\psi_j\in B_{2}} \gamma_{j}t_j^kE_{1j}\in T\mathcal{F}^k(j^k\psi)\subseteq V\subseteq T\mathcal{G}^k(j^k\psi)$ be a nonzero vector and denote $J'=\{j;\ \gamma_{j}\neq 0\}\subseteq I$. Remark that $l\not\in J'$ for all $l\in L_k$.

Choosing $j_0\in J'=:J'_{j_0}$, taking the element $w_l:=v_l-\frac{\alpha_{lj_0}}{\gamma_{j_0}}w \in T\mathcal{G}^k(j^k\psi)$ and setting $J'_l=\{j\in J_l\setminus\{j_0\};\ \alpha_{lj}\gamma_{j_0}\neq\alpha_{lj_0}\gamma_{j}\}$ for all $l\in L_k$, it follows by Corollary \ref{CorDifTg} that $F_{J'_l}(\underline{k})\neq \emptyset$. In addition for each $l\in L':=L_k\cup\{j_0\}$ we have that $l \in J'_l \setminus \bigcup_{i\in L'\setminus\{l\}}J'_i$ and $F_{J'_l}(\underline{k})\neq \emptyset$. In this way, $L'\supsetneq L_k$ satisfies the condition $(\star)$ in Proposition \ref{PropElim}, which contradicts the maximality of the set $L_k$.

Hence, $\psi=\psi'$ and the $\mathcal{G}$-normal form of $\varphi$,  as (\ref{normalBlocks}), is unique.
\cqd

\begin{remark}
Notice that the parameters $a_{ij}$ in the $\mathcal{G}$-normal form (\ref{normalBlocks}) are not necessarily free because, in order to guarantee that the multigerm is in $\Sigma_{\Gamma,\Lambda_{\mathcal{G}}}$, some algebraic conditions must be imposed on the coefficients. The conditions imposed to have the semiring $\Gamma$ are easily determined by the nonzero coefficients of the terms with characteristic exponents and some conditions to get the intersection multiplicity of each pair of branches. On the other hand, the algebraic conditions related with $\Lambda_{\mathcal{G}}$ can be obtained applying the algorithm presented in \cite{carvalho2}.
\end{remark}

As we remarked in Section 2, the analytic equivalence of plane curves with $r$ irreducible components can be translated by the $\mathcal{S}\times\mathcal{A}$-action on $\mathcal{P}$. Any element in $\mathcal{P}$ is $\mathcal{S}\times\mathcal{A}$-equivalent to a multigerm given in a Puiseux block form and the $\mathcal{A}$-action on such elements splits according to the action of subgroups given in Lemma \ref{subgroups}, namely $\mathcal{H}\circ\tilde{\mathcal{A}_1}$, $\mathcal{H}\circ\mathcal{A}_1$ or $\mathcal{H}'\circ\mathcal{A}_1$.

In Theorem \ref{A1normal-form} we presented $\mathcal{G}$-normal forms for Puiseux block form multigerms with $\mathcal{G}=\tilde{\mathcal{A}_1}$ or
$\mathcal{G}=\mathcal{A}_1$. In this way, to present $\mathcal{A}$-normal forms of such multigerms we have to consider the action of $\mathcal{H}$ or $\mathcal{H}'$, that we call {\it homothety actions}, on the respective $\mathcal{G}$-normal forms. Recall that an action by a subgroup of $\mathcal{H}$ does not introduce nor eliminate terms in a multigerm.

Let us consider $\varphi=[\varphi_1,\ldots ,\varphi_r]\in\mathcal{P}$ given as (\ref{normalBlocks}) and the parameter vector
\begin{equation}\label{param-vector}
p=(a_{1,n_1+1}, \ldots, a_{1,k_0},a_{2,n_2+1}, \ldots, a_{2,k_0},\ldots ,a_{r,n_r+1}, \ldots, a_{r,k_0}),
\end{equation}
where $k_0=\min\{z\in\mathbb{N};\ \underline{z}\in\varrho+\mathbb{N}^r\}$ (see Remark \ref{propriedade-max}).

Given an element $(\rho_1,\ldots ,\rho_r,\sigma)\in\mathcal{H}$  where $\rho_i(t_i)=u_it_i$ and $\sigma(X,Y)=(\alpha X,\gamma Y)$ with $\alpha,\gamma,u_i\in\mathbb{C}^*$ for $i\in I$, in order to keep  $\psi=\sigma\circ\varphi\circ\rho^{-1}$ as in (\ref{normalBlocks}), that is, $j^{n_i}(\varphi_i)=j^{n_i}(\psi_i)$ for all $i\in I$, we must have $$u_i=\alpha^{-\frac{1}{n_i}}\ \ \mbox{if}\ \ \varphi_i\not\in B_2,\ \ \ \ \ \ \ \ u_i=\gamma^{-\frac{1}{n_i}}\ \ \mbox{if}\ \ \varphi_i\in B_2$$ and $\theta_l\gamma\alpha^{-1}=\theta_l$ for all $\varphi_i\in B_l$ with $l\neq 2$. Remark that $\alpha=\gamma$ for $(\rho_1,\ldots ,\rho_r,\sigma)\in\mathcal{H}'\subset\mathcal{H}$, that is, if we have $s>2$ blocks.

The corresponding parameter vector of $\psi$ is
$
(b_{1,n_1+1}, \ldots, b_{1,k_0},\ldots ,b_{r,n_r+1}, \ldots, b_{r,k_0})$ with
\begin{equation}\label{coef}
b_{ij}=\left \{ \begin{array}{ccc}\gamma\alpha^{-\frac{j}{n_i}}a_{ij} & \mbox{if} & \psi_i\in B_1 \\ \alpha\gamma^{-\frac{j}{n_i}}a_{ij} & \mbox{if} & \psi_i\in B_2\end{array}\right .\ \mbox{for}\ s\leq 2,\ \ \mbox{and} \ \ b_{ij}=\alpha^{\frac{n_i-j}{n_i}}a_{ij}\ \ \mbox{for all}\ \ i\in I,\  \mbox{if}\ s\geq 3.\end{equation}

Notice that the parameter vector $p$ is null if and only if
\begin{equation}\label{basic}
\varphi_1=\psi_1=(t_1,0),\ \ \ \varphi_2=\psi_2=(0,t_2),\ \ \ \varphi_i=\psi_i=(t_i,\theta_it_i)\ \ \mbox{for all}\ \ 3\leq i\leq r,
\end{equation}
with $\theta_i\neq\theta_j$ for $i\neq j$.

On the other hand, if $p$ is a nonzero vector and $a_{i_0j_0}$ is the first nonzero  coordinate, then $\psi_i$ (and $\varphi_i$) for all $i<i_0$ is given as (\ref{basic}) and we can choose $\alpha\in\mathbb{C}^*$ in (\ref{coef}) corresponding to $a_{i_0j_0}$ in order to normalize the corresponding coefficient in $\psi$, that is, to obtain $b_{i_0j_0}=1$.

Eventually the equations (\ref{coef}) allow us to normalize one more coefficient, but for this purpose a finer analysis would be necessary. However, the above description it is sufficient to present the following theorem that gives a solution (up to $\mathcal{S}$-action) for the analytical equivalence problem for plane curves with $r$ irreducible components.

\begin{theorem}($\mathcal{A}$-{\sc normal form})\label{Anormal-form}
    Let $\varphi=[\varphi_1,\ldots, \varphi_r]\in \Sigma_{\Gamma,\Lambda_{\mathcal{G}}}$ be a multigerm associated to the curve $\mathcal{C}$. Then $\varphi$ is $\mathcal{A}$-equivalent to
    $\psi=[\psi_1,\ldots, \psi_r]$ as (\ref{normalBlocks}) with parameter vector $p=(b_{1,n_1+1}, \ldots, b_{r,k_0})$ such that

    If $p=(0,\ldots ,0)$ then $\psi_1=(t_1,0),\ \psi_2=(0,t_2),\  \psi_i=(t_i,\theta_it_i)\ \mbox{for all}\ 3\leq i\leq r$.

    If $p\neq (0,\ldots ,0)$ then $b_{i_0j_0}=1$ is the first nonzero coordinate of $p$. In addition, $\psi'\in \Sigma_{\Gamma,\Lambda_{\mathcal{G}}}$ as (\ref{normalBlocks}) with parameter vector $p'=(b'_{1,n_1+1}, \ldots, b'_{r,k_0})$ and first nonzero coordinate $b'_{i_0j_0}=1$ is $\mathcal{A}$-equivalent to $\psi$ if and only if we have the following conditions according to the number $s$ of blocks:
    \begin{center}
        \begin{tabular}{|c|l|}
            \hline
            & \\
             $s=1$ & There exists $c\in\mathbb{C}^*$ such that $b_{ij}=c^{\frac{jn_{i_0}-j_0n_i}{n_in_{i_0}}}b'_{ij}$ for all coordinates of $p$ and $p'$. \\ \hline

             & If $\psi_{i_0}\in B_1$ there exists $c\in\mathbb{C}^*$ such that:\\

        &    \hspace{1.5cm} $b_{ij}=c^{\frac{jn_{i_0}-j_0n_i}{n_in_{i_0}}}b'_{ij}$ for all coordinates of $p$ and $p'$ with $\psi_i\in B_1$; \\

     $s=2$  &    \hspace{1.5cm} $b_{ij}=c^{\frac{jj_0-n_{i_0}n_i}{n_in_{i_0}}}b'_{ij}$ for all coordinates of $p$ and $p'$ with $\psi_i\in B_2$. \\ \cline{2-2}
        & \\
        &    If $\psi_{i_0}\in B_2$ there exists $c\in\mathbb{C}^*$ and  $b_{ij}=c^{\frac{jn_{i_0}-j_0n_i}{n_in_{i_0}}}b'_{ij}$ for all entries of $p$ and $p'$. \\ \hline
        & \\
     $s\geq 3$ & There exists $c\in\mathbb{C}^*$; $c^{\frac{j_0-n_{i_0}}{n_{i_0}}}=1$ and $b_{ij}=c^{\frac{j-n_{i}}{n_i}}b'_{ij}$ for all entries of $p$ and $p'$. \\
            \hline
        \end{tabular}
    \end{center}
\end{theorem}
\Dem By Theorem \ref{A1normal-form} we obtain a unique $\mathcal{G}$-normal form for $\varphi$ as (\ref{normalBlocks}) and, as above described, if the parameter vector is nonzero then the $\mathcal{G}$-normal form is $\mathcal{H}$-equivalent (or $\mathcal{H}'$-equivalent if $s\geq 3$) to $\psi=[\psi_1,\ldots, \psi_r]$ as (\ref{normalBlocks}) with parameter vector $p=(b_{1,n_1+1}, \ldots, b_{r,k_0})$ such that its first nonzero coordinate is $b_{i_0j_0}=1$.

Now it is sufficient to describe the parameter vector $p'=(b'_{1,n_1+1}, \ldots, b'_{r,k_0})$ of a multigerm  $\psi'\in \Sigma_{\Gamma,\Lambda_{\mathcal{G}}}$ as (\ref{normalBlocks}) with the first nonzero coordinate $b'_{i_0j_0}=1$ in the orbit $\mathcal{H}(\psi)$ if $s\leq 2$ or in $\mathcal{H}'(\psi)$ if $s\geq 3$. This will be done using the equations (\ref{coef}).

If $s=1$, by (\ref{coef}) we have $\gamma \alpha^{-\frac{j_0}{n_{i_0}}}=1$ and $\gamma \alpha^{-\frac{j}{n_{i}}}b_{ij}=b'_{ij}$. Consequently, $\gamma=\alpha^{\frac{j_0}{n_{i_0}}}$ and $b_{ij}=\alpha^{\frac{jn_{i_0}-j_0n_i}{n_in_{i_0}}}b'_{ij}$ for all coordinates of $p$ and $p'$.

For multigerms with two blocks it follows, by (\ref{coef}), the equations $b'_{ij}=\gamma \alpha^{-\frac{j}{n_{i}}}b_{ij}$ if $\psi_i\in B_1$ and $b'_{ij}=\alpha\gamma^{-\frac{j}{n_{i}}}b_{ij}$ if $\psi_i\in B_2$.
If $\psi_{i_0}\in B_1$, then $\gamma=\alpha^{\frac{j_0}{n_{i_0}}}$ and $b_{ij}=\alpha^{\frac{jn_{i_0}-j_0n_i}{n_in_{i_0}}}b'_{ij}$ for coordinates of $p$ and $p'$ with $\psi_i\in B_1$ and $b_{ij}=\alpha^{\frac{jj_0-n_{i_0}n_i}{n_in_{i_0}}}b'_{ij}$ for all coordinates of $p$ and $p'$ with $\psi_i\in B_2$.
If $\psi_{i_0}\in B_2$, then $\alpha=\gamma^{\frac{j_0}{n_{i_0}}}$ and $b_{ij}=\gamma^{\frac{jn_{i_0}-j_0n_i}{n_in_{i_0}}}b'_{ij}$ for all entries of $p$ and $p'$.

Finally, for $s\geq 3$ we have $\alpha^{\frac{n_{i_0}-j_0}{n_{i_0}}}=1$ and $\alpha^{\frac{n_i-j}{n_i}}b_{ij}=b'_{ij}$ for all $i\in I$ and $n_i+1\leq j\leq k_0$.
\cqd

The above theorem provides a method to decide if two plane curves are or not analytically equivalent. In fact, considering an associated multigerm to a plane curve $\mathcal{C}$ with $r$ irreducible components we compute their $\mathcal{A}$-invariants $\Gamma$ and $\Lambda_{\mathcal{G}}$ using the appropriate group $\mathcal{G}$ and the algorithms in \cite{carvalho2} for instance. Applying the results in this section we obtain an $\mathcal{A}$-normal form $\varphi=[\varphi_1,\ldots ,\varphi_r]$ for the multigerm of $\mathcal{C}$.

Given another plane curve $\mathcal{C}'$ with multigerm associated $\psi=[\psi_1,\ldots ,\psi_r]$, if the $\mathcal{A}$-invariants $\Gamma^{\pi}$ and $\Lambda^{\pi}_{\mathcal{G}}$ of $\pi(\psi)$
are distinct to the corresponding ones of $\varphi$ for every $\pi\in\mathcal{S}$, then $\mathcal{C}$ and $\mathcal{C}'$ are not analytically equivalent. If it is not the case, we consider the subgroup $\mathcal{R}$ of $\mathcal{S}$ such $\Gamma^{\pi}=\Gamma$ and $\Lambda^{\pi}_{\mathcal{G}}=\Lambda_{\mathcal{G}}$ for all $\pi\in\mathcal{R}$. The curves $\mathcal{C}$ and $\mathcal{C}'$ are analytically equivalent if and only if there exists $\pi\in\mathcal{R}$ such that $\varphi$ and an $\mathcal{A}$-normal form for $\pi(\psi)$ satisfy the corresponding condition in Theorem \ref{Anormal-form}.

\section{Particular cases}

In this final section we will recover some known results regarding the analytic classification of plane curves using Theorems \ref{A1normal-form} and  \ref{Anormal-form} and we illustrate how we can apply them in particular topological classes.

\subsection{Irreducible case and bigerms}

For the irreducible case, i.e., $s=r=1$, we must consider $\mathcal{G}=\tilde{\mathcal{A}}_1$ and the normal form described in Theorem \ref{A1normal-form} is rewritten as
$\varphi=(t^n,\sum_{j>n\atop L_j=\emptyset}c_jt^j).$

Notice that $L_j=\emptyset$ if and only if $F_{\{1\}}(j)=\emptyset$, that is, $j\not\in\Lambda_{\mathcal{G}}$.

If $n=1$, then $\Lambda_{\mathcal{G}}=\overline{\mathbb{N}}\setminus\{0,1\}$ and $\varphi$ is $\mathcal{A}$-equivalent to $(t,0)$.

Let us suppose that $n\geq 2$ and $\Gamma\cap\mathbb{N} =\langle n,v_1,\ldots ,v_g\rangle$. Thus,
$\varphi=(t^n,c_{v_1}t^{v_1}+\sum_{j\not\in \Lambda_{{\mathcal{G}}}}c_jt^j)$ with $c_{v_1}\neq 0$ and, by Theorem \ref{Anormal-form}, we have $\varphi\underset{\mathcal{A}}{\sim}(t^n,t^{v_1}+
\sum_{j\not\in \Lambda_{{\mathcal{G}}}}b_jt^j)$.

If $j\in\Lambda_{\mathcal{G}}$ for every $j>v_1$, then $\varphi\underset{\mathcal{A}}{\sim}(t^n,t^{v_1})$. By \cite{torcao}, this is equivalent to claim that $\Lambda\setminus\Gamma=\emptyset$.

If there exists $j\not\in\Lambda_{\mathcal{G}}$ for some $j>v_1$ then setting $\lambda=\min\{j>v_1;\ j\not\in\Lambda_{\mathcal{G}}\}$ we have \begin{equation}
\label{IrredNormalForm}
\varphi \ \underset{\mathcal{A}}{\sim}\ (t^n,t^{v_1}+b_{\lambda}t^{\lambda}+
\sum_{\lambda<j\not\in \Lambda_{{\mathcal{G}}}}b_jt^j)
\end{equation} with $b_{\lambda}\neq 0$ and $\lambda+n=\nu(\varphi^*(\omega_0))$ where $\omega_0:=v_1YdX-nXdY$.

\begin{remark}
    Notice that $\{j>\lambda;\ j\not\in\Lambda_{\mathcal{G}}\}=\{j>\lambda;\ j\not\in\Lambda-n\}$. In fact, by definition we have $\Lambda_{\mathcal{G}}\subset\Lambda -n$. On the other hand, given $\omega=(\alpha_0+\alpha_1X+\alpha_2Y+h_1)dX+(\alpha_3+\alpha_4X+h_2)dY\in\Omega^1$ with $\alpha_i\in\mathbb{C}$ and $\omega_1=h_1dX+h_2dY\in\Omega_{\mathcal{G}}$ in order to have $\nu(\varphi^*(\omega))-n>\lambda$ we must have $\alpha_0=\alpha_1=\alpha_3=0$ and $n\alpha_2+v_1\alpha_4=0$. In this way, $\omega=\frac{\alpha_2}{v_1}\omega_0+\omega_1$, as  $\nu(\varphi^*(\omega_0))-n=\lambda\neq \nu(\varphi^*(\omega_1))-n\in\Lambda_{\mathcal{G}}$ (by (\ref{IrredNormalForm})), we have $\nu(\varphi^*(\omega))-n=\nu(\varphi^*(\omega_1))-n\in\Lambda_{\mathcal{G}}$.
\end{remark}

By the above remark, if $\Lambda\setminus\Gamma\neq\emptyset$, that is, there exists $\lambda=\min\{j>v_1;\ j\not\in\Lambda_{\mathcal{G}}\}$ then $\varphi\underset{\mathcal{A}}{\sim}(t^n,t^{v_1}+b_{\lambda}t^{\lambda}+
\sum_{\lambda<j\not\in \Lambda-n}b_jt^j)$ with $b_{\lambda}\neq 0$.

Taking $c\in\mathbb{C}^*$ with $c^{\lambda-v_1}=b_{\lambda}$, $\rho(t)=ct$ and $\sigma(X,Y)=(c^nX,c^{v_1}Y)$ we obtain $(\rho,\sigma)\cdot\varphi=(t^n,t^{v_1}+t^{\lambda}+
\sum_{\lambda<j\not\in \Lambda-n}a_jt^j)$,
that is, the normal form presented in Theorem 2.1 in \cite{HefezHern}.

For bigerms, we have two distinct cases according to the number $s\in\{1,2\}$ of Puiseux blocks.

By Theorem \ref{A1normal-form} we get
$$\varphi_1=\left ( t_1^{n_1},\sum_{j>n_1\atop 1\not\in L_j}a_{1j}t_1^j\right )\ \ \ \ \mbox{and}\ \ \ \ \varphi_2=\left \{
\begin{array}{ll}
\left (  t_2^{n_2},\sum_{j>n_2\atop 2\not\in L_j}a_{2j}t_2^j\right ) & \mbox{if}\ s=1;\\ & \\
\left ( \sum_{j>n_2\atop 2\not\in L_j}a_{2j}t_2^j, t_2^{n_2}\right )  & \mbox{if}\ s=2.
\end{array}\right .$$

Notice that $1\not\in L_j$ if and only if $F_{\{1\}}(\underline{j})=F_{\{1,2\}}(\underline{j})=\emptyset$, that is, $F_{\{1\}}(\underline{j})=\emptyset$ and $\underline{j}$ is not a maximal point of $\Lambda_{\mathcal{G}}$. On the other hand $2\not\in L_j$ means that $F_{\{2\}}(\underline{j})=\emptyset$.

Hence, for curves with two branches, Theorem \ref{A1normal-form} recovers Theorem 6 in \cite{HefezHernRod}.

\subsection{Singular ordinary point}

Let us consider a plane curve $\mathcal{C}$ with a singular ordinary point at
the origin of multiplicity $r>1$, that is, $\mathcal{C}$ has $r$
irreducible regular components given by $f_i\in\mathcal{M}\setminus
\mathcal{M}^2$ such that the slopes $\theta_i$ of their tangent lines satisfy
$\theta_i\neq\theta_j$ for all $i,j\in I,\ i\neq j$.

Applying Proposition 3.9 of \cite{carvalho1} for this case we obtain that the value semiring $\Gamma$ is minimally generated by $$v_i=\nu(f_i)=(1,\ldots ,1,\infty ,1,\ldots ,1);\ i\in I$$ where the symbol ``$\infty$'' occupies the $i^{th}$-coordinate. In particular, any element $\gamma\in\Gamma$ can be written as (\ref{semiring}), i.e.,
$\gamma=\inf\left \{\sum_{i\in I}a_{1i}v_i,\ldots,\sum_{i\in I}a_{ri}v_i\right \}.$ If $a_{li}\neq 0$ for all $i\in I$, then $\sum_{i\in I}a_{li}v_i=(\infty,\ldots ,\infty)$ which is irrelevant for the computation of $\gamma$. So, we can assume $a_{li}=0$ for every $l\in I$ and some $i\in I$ which depends on $l$.

Notice that $\mu_i=0$ and $\nu_j(f_i)=1$ for all $i,j\in I$ with $i\neq j$, then the conductor of $\Gamma$ is $\kappa=(r-1,\ldots ,r-1)$. In particular, $F^{\Gamma}_{J}(\underline{k})\neq\emptyset$ for every $k\geq r-1$ and all $\emptyset\neq J\subseteq I$.

By Propositions \ref{trunc} and \ref{di}, for $r\in\{1,2,3\}$ the plane curve $\mathcal{C}$ admits a multigerm $\varphi$ that is $\mathcal{A}$-equivalent to
$[(t_1,0)], [(t_1,0),(0,t_2)]$ or $[(t_1,0),(0,t_2),(t_3,t_3)]$, respectively.

In what follows we consider $r\geq 4$. Notice that the appropriated subgroup to apply  Theorem \ref{A1normal-form} is $\mathcal{G}=\mathcal{A}_1$.

Using the generators of the value semiring $\Gamma$ we can describe the non-empty fibers $F^{\Gamma}_J(\underline{k})$ for $0\leq k<r$.

\begin{lemma}
    For a singular ordinary point with multiplicity $r$ the $J$-fiber of $\underline{k}\in\mathbb{N}^r$ of $\Gamma$ is non-empty if and only if $\sharp J\geq r-k$. Moreover, for any $J\subseteq I$ with $\sharp J\geq r-k$ we have $$F_{J}^{\Gamma}(\underline{k})=\{\underline{k}+\alpha;\ \ \alpha=(\alpha_1,\ldots ,\alpha_r)\in\overline{\mathbb{N}}^r\ \ \mbox{with}\ \ \alpha_j=0\ \ \mbox{if}\ \ j\in J\ \ \mbox{and}\ \ \alpha_j\neq 0\ \ \mbox{if}\ \ j\in I\setminus J\}.$$
\end{lemma}
\Dem  As the conductor of $\Gamma$ is $\kappa=(r-1,\ldots ,r-1)$ it is sufficient
    to describe the fiber $F^{\Gamma}_{J}(\underline{k})$ for every
    $\underline{k}=(k,\ldots ,k)\in\mathbb{N}^r,\ 0\leq k\leq r-1$ and
    every $J\subseteq I$ with $J\neq\emptyset$.

    Firstly, notice that if $F^{\Gamma}_{J}(\underline{k})\neq \emptyset$ then $\sharp J\geq r-k$. In fact, let $\gamma\in F^{\Gamma}_{J}(\underline{k})$. Since $\gamma \in \Gamma$, there exists $h\in\mathbb{C}\{X,Y\}$, namely of multiplicity $n_h$, such that $\gamma=\nu(h)$. If $\sharp J< r-k$, as $n_i=1$, $\theta_i\neq\theta_l$ for all $i,l\in I$ with $i\neq l$, we must have $m_h\leq n_h\leq\nu_j(h)=k<r-\sharp J$, for all $j \in J$, where $m_h$ is the number of distinct  tangent lines of $h$. But in this way, there exists $i\in I\setminus J$ such that $f_i$ and $h$ do not share the same tangent line. So, $\nu_i(h)=n_h\leq k$ and therefore $\gamma\not\in F^{\Gamma}_{J}(\underline{k})$ which is a contradiction.

    Now we will describe $F^{\Gamma}_{J}(\underline{k})$ with $\sharp J\geq r-k$.

    For $0\leq k\leq r-1$ we take any $J_0\subseteq I$
    with $\sharp J_0=r-k\ (\geq 1)$ and let us consider the sets $I\setminus J_0=\{i_1,\ldots ,i_k\}$ and
    $\{i_1,\ldots, \widehat{i_m},\ldots ,i_k\}:=\{i_1,\ldots
    ,i_k\}\setminus\{i_m\}$.

    Fixing $i_0\in J_0$ we set
    $$\zeta_0=\sum_{l\in \{i_1,\ldots ,i_k\}}\nu(f_l),\ \ \ \ \zeta_m=\sum_{l\in \{i_1,\ldots ,\widehat{i_m},\ldots i_k\}}\nu(f_l)+(1+\alpha_{m})\nu(f_{i_0})\in\overline{\mathbb{N}}^r,$$ with $1\leq m\leq k$ and $\alpha_1,\ldots ,\alpha_k\in\mathbb{N}$. Notice that the coordinates of $\zeta_0$ and $\zeta_m$, respectively, satisfy
    \begin{center}
        \begin{tabular}{ccc}
            $\zeta_{0j}=\left \{\begin{array}{lcl} k & \mbox{if} & j\in J_0 \\
            \infty & \mbox{if} & j\in I\setminus J_0,\end{array} \right .$
            & & $\zeta_{mj}=\left \{\begin{array}{lcl}
            \infty & \mbox{if} & j\in \{i_0,i_1,\ldots ,\widehat{i_m},\ldots i_k\} \\
            k+\alpha_{m} & \mbox{if} &  j\not\in \{i_0,i_1,\ldots ,\widehat{i_m},\ldots i_k\}.\end{array} \right .$ \end{tabular}
    \end{center}

    So, $\gamma=(\gamma_1,\ldots ,\gamma_r):=\inf\left \{\zeta_0,\zeta_1,\ldots ,\zeta_k\right
    \}\in\Gamma$ is such that $\gamma_j=k$
    if $j\in J_0$, $\gamma_j=k+\alpha_j$ with $\alpha_j\in\mathbb{N}$ for
    $j\in I\setminus J_0$.

    Hence, the above construction ensures that for any $J\subseteq I$ with $\sharp J\geq r-k$ we have
    $$F_{J}^{\Gamma}(\underline{k})=\{\underline{k}+\alpha;\ \ \alpha=(\alpha_1,\ldots ,\alpha_r)\in\overline{\mathbb{N}}^r\ \ \mbox{with}\ \ \alpha_j=0\ \ \mbox{if}\ \ j\in J\ \ \mbox{and}\ \ \alpha_j\neq 0\ \ \mbox{if}\ \ j\in I\setminus J\}.$$
\cqd

Remark that $\nu(\mathcal{M}^2)-\underline{1}\subseteq \Lambda_{\mathcal{A}_1}$ which implies that its corresponding fibers satisfy $\emptyset\neq F_J^{\Gamma}(\underline{k+1})-\underline{1}\subseteq F_J(\underline{k})$ for all $k\geq 2$ and  $\sharp J \geq r-(k+1)$.

Taking $L=\{1,\ldots ,\min\{k+2,r\}\}$ and $J_i=\{i,k+3,\ldots ,r\}$ for $i\in L$ in Proposition \ref{PropElim}, then Corollary \ref{CorElim} allows us to eliminate the term corresponding to $t_i^k$ of the component $\varphi_i$ for all $i=1, \dots, \min\{k+2,r\}$ preserving the analytic class of $\mathcal{C}$.
With this approach we obtain Theorem 1.1 in \cite{genzmerpaul1}:

\begin{proposition}
    Any curve with a singular ordinary point with $r$ branches is analytically equivalent to a curve defined by a multigerm
    $\psi=[\psi_1, \ldots, \psi_r]$ with
    \begin{equation}
    \label{CurvRegTransv}
    \psi_1=(t_1,0),\ \ \psi_2=(0,t_2), \ \  \psi_i=\left(t_i, \ \theta_{i}t_i +
    \displaystyle\sum_{j=2}^{i-3}a_{ij}t_i^j\right),\ i=3, \dots,
    r,
    \end{equation}
    where $\theta_3=1$ and $\theta_i\neq\theta_l$ for $i\neq l$.
\end{proposition}

For $r\geq 4$, Granger in \cite{Granger} considers the generic parameter vector to obtain $\mathcal{A}$-normal forms and to compute the dimension of the moduli space. The same result was obtained with
other methods by Genzmer and Paul in \cite{genzmerpaul1}. From now on, we consider the generality hypothesis and we will apply our techniques with the purpose to compute the dimension of the moduli space.

The generality of the coefficients give us a symmetry on the coordinates of $\Lambda_{\mathcal{A}_1}$. In this way, for $k\geq 2$ the set $L_k\subseteq I$ in Theorem \ref{A1normal-form} is $L_k=\{1,2,\ldots ,e(k)\}$ for some $4\leq e(k)\leq r$. By Remark \ref{lk-diff}, $e(k)$ is the maximum integer such that for any $b_{ik}\in\mathbb{C}$ with $1\leq i\leq e(k)$ there exists ${\rm w}\in \mathcal{I}_{\mathcal{A}_1}$ with $j^k{\rm w}=(b_{1k}t_1^{k},\ldots ,b_{rk}t_r^{k})$.

Recall that ${\rm w}=\left (t_1^{-1}\varphi^*_1(\omega),\ldots ,t_r^{-1}\varphi^*_r(\omega)\right )\in \mathcal{I}_{\mathcal{A}_1}$ and $t_i^{-1}\varphi^*_i(\omega)=\varphi^*_i(p)\cdot x'_i+\varphi^*_i(q)\cdot y'_i\in\mathbb{C}\{t_i\}$ for $i\in I$ where $\omega=pdX+qdY\in\Omega^1$ and $p,q\in\langle X,Y\rangle^2$.

In what follows we denote
\begin{equation}
\label{omegah}
\omega_h=\left(\displaystyle\sum_{m=0}^{h}\alpha_{h,m}X^{h-m}Y^m\right)dX+\left(\displaystyle
\sum_{m=0}^{h}\beta_{h,m}X^{h-m}Y^m\right)dY \in \Omega^1
\end{equation} with $\alpha_{h,m},\beta_{h,m}\in\mathbb{C}$. Notice that given $\varphi=[\varphi_1,\ldots ,\varphi_r]$ as (\ref{CurvRegTransv}) and $\omega \in \Omega^1$, if we want to obtain $\nu_i(\omega)=k$ and consequently to determine $e(k)$, it is sufficient to consider $w=\sum_{h=2}^{k}w_h$.

\begin{lemma}
    \label{LemaCurvRegTransv}
With the above notations, if $r\geq 4$ we have $e(2)=4$ and $e(k)=\min\{2k+1,r\}$ for every $k\geq 3$.
\end{lemma}
\Dem Considering a generic multigerm $\varphi$ given as in (\ref{CurvRegTransv}) if $\omega=\omega_2$, that is $k=2$, then
$$j^2(t_1^{-1}\varphi_1^{*}(\omega))=\alpha_{2,0}t_1^2,\ \ j^2(t_2^{-1}\varphi_2^{*}(\omega))=\beta_{2,2}t_2^2,$$
$$j^2(t_i^{-1}\varphi_i^{*}(\omega))=\left( \alpha_{2,0}+(\alpha_{2,1}+\beta_{2,0})\theta_i+(\alpha_{2,2}+\beta_{2,1})\theta_i^2+\beta_{2,2}\theta_i^3\right )t_i^2,\ \ \mbox{for all}\ \ 3\leq i\leq r.$$
Now, it is immediate to verify that the maximum integer $e(2)$ such that the system $j^2(t_i^{-1}\varphi_i^{*}(\omega))=b_{i2}t_i^2$ admits solution for any $b_{i2}\in\mathbb{C}$ with $1\leq i\leq e(2)$ is precisely $e(2)=4$.

We will show that for any $3\leq l\leq \left [ \frac{r-1}{2}\right ]$ if $j^{l}(t_i^{-1}\varphi_i^{*}(\omega_{l-1}+\omega_{l}))=0$ for every $i\in I$, then $w_{l-1}=0$.

Notice that for $w_{h}$ as (\ref{omegah}) with $h\geq 2$, we get
{\small \begin{equation}
    \label{Omegak}
    \begin{array}{l}
    j^{h+1}(t_1^{-1}\varphi_1^{*}(\omega_h))=\alpha_{h,0}t_1^h,\\ j^{h+1}(t_2^{-1}\varphi_2^{*}(\omega_h))=\beta_{h,h}t_2^h, \\
     j^{h+1}(t_i^{-1}\varphi_i^{*}(\omega_h))=\left(\alpha_{h,0}+\sum_{m=1}^h(\alpha_{h,m}+\beta_{h,m-1})\theta_{i}^m +\beta_{h,h}\theta_i^{h+1}\right )t_i^h+ \\ \hspace{3cm} +\left ( \sum_{m=0}^{h-1}((m+1)\alpha_{h,m+1}+(m+2)\beta_{h,m})\theta_{i}^{m}+(h+2)\beta_{h,h}\theta_i^{h}\right )a_{i2}t_i^{h+1}\\ \hspace{12cm} \mbox{for}\ \ 3\leq i\leq r.
    \end{array}
    \end{equation}
}

The condition $j^{l}(t_i^{-1}\varphi_i^{*}(\omega_{l-1}+\omega_{l}))=0$ for every $i\in I$ implies that $\alpha_{l-1,0}=\alpha_{l,0}=\beta_{l-1,l-1}=\beta_{l,l}=0$ and the coefficient of $t^{l-1}$ in $t_i^{-1}\varphi_i^{*}(\omega_{l-1}+\omega_{l})$ vanishes for every $3\leq i\leq r$, that is,
$$
\left(\begin{array}{cccc}
\theta_{3} & \theta_{3}^2 & \cdots & \theta_{3}^{l}  \\
\vdots & \vdots &  & \vdots  \\
\theta_{r} & \theta_{r}^2 & \cdots & \theta_{r}^{l}  \\
\end{array}\right) \cdot \left(\begin{array}{c}
\alpha_{l-1,1}+\beta_{l-1,0}\\
\vdots \\
\alpha_{l-1,l-1}+\beta_{l-1,l-2}\\
\end{array}\right)=\left(\begin{array}{c}
0\\
\vdots \\
0\\
\end{array}\right).
$$
So, $\alpha_{l-1,m}=-\beta_{l-1,m-1}$ for $1\leq m\leq l-1$. In addition, as the coefficient of $t^{l}$ in $t_i^{-1}\varphi_i^{*}(\omega_{l-1}+\omega_{l})$ is null, by (\ref{Omegak}) we must have
\begin{equation}
\label{SistRegTrans}
\left(\begin{array}{ccccccc}
a_{32} & \theta_{3}a_{32} & \cdots & \theta_{3}^{l-2}a_{32} & \theta_{3} & \cdots & \theta_{3}^{l}  \\
\vdots & \vdots &  & \vdots & \vdots & & \vdots \\
a_{r2} & \theta_{3}a_{32} & \cdots & \theta_{r}^{l-2}a_{r2} & \theta_{r} & \cdots & \theta_{r}^{l}  \\
\end{array}\right) \cdot \left(\begin{array}{c}
\beta_{l-1,0}\\
\vdots \\
\beta_{l-1,l-2}\\
\alpha_{l,1}+\beta_{l,0}\\
\vdots \\
\alpha_{l,l}+\beta_{l,l-1}
\end{array}\right)=\left(\begin{array}{c}
0\\
\vdots \\
0\\
\end{array}\right).
\end{equation}
As we suppose $\varphi$ generic, the coefficient matrix has maximal rank $\min\{2l-1,r-2\}$ and we get the solution $\beta_{l-1,m}=0$ for $0\leq m\leq l-2$ and $\alpha_{l,m}=-\beta_{l,m-1}$ for $1\leq m\leq l$. Consequently, $w_{l-1}=0$.

In order to obtain $j^{k}(t_i^{-1}\varphi_i^{*}(\sum_{h=2}^{k}\omega_h))=b_{ik}t_i^{k}$ it is necessary that $j^{l}(t_i^{-1}\varphi_i^{*}(\omega_h))=0$ for each $h=2,\ldots, k-1$ and for all $l=h,\ldots, k-1$. From this and the above discussion, it follows by induction that $\omega_h=0$ for $2\leq h\leq k-2$. So, it is sufficient to consider $\omega=\omega_{k-1}+\omega_{k}$.

By a similar computation as in (\ref{Omegak}) using $l=k$ we get the system
$$\left \{\begin{array}{l}
\alpha_{k,0}=b_{1k}\\
\beta_{k,k}=b_{2k}\\
\alpha_{k,0}+\sum_{m=0}^{k-2}\beta_{k-1,m}\theta_i^{m}a_{i2}+\sum_{m=1}^{k}(\alpha_{k,m}+\beta_{k,m-1})\theta_i^m+\beta_{k,k}\theta_i^{k+1}=b_{ik}\ \ \ \mbox{for}\ 3\leq i\leq r.
\end{array}\right .
$$

Substituting the two first solutions in the other equations, we obtain a system with same coefficient matrix as (\ref{SistRegTrans}). So, generically we can solve at most $e(k)=2+\min\{2k-1,r-2\}=\min\{2k+1,r\}$ such equations.
\cqd

As a consequence we recover the analytical normal form for regular transversal branches in the generic case, as given in \cite{Granger} (Proposition 2) and \cite{genzmerpaul1}.

\begin{proposition}
Generically any plane curve with a singular ordinary point and $r \geq 4$ is analytically equivalent to a curve defined by a multigerm
    $\psi=[\psi_1, \ldots, \psi_r]$ with
$$\begin{array}{c}
\psi_1=(t_1,0),\ \ \psi_2=(0,t_2), \ \  \psi_3=(t_3,t_3),\ \ \psi_4=(t_4,\theta_4t_4)\vspace{0.2cm}\\
\psi_i=\left(t_i, \ \theta_{i}t_i + a_{i2}t_i^2+
\displaystyle\sum_{k=3}^{\left [ \frac{i-2}{2}\right  ]}a_{ik}t_i^k\right),\ 5\leq i\leq r.
\end{array}
$$
\end{proposition}
\Dem The previous lemma gives us $L_2=\{1,2,3,4\}$ and $L_k=\{1,\ldots ,\min\{2k+1,r\}\}$ for $k\geq 3$. So, by Theorem \ref{A1normal-form}, we have the above normal form for $1\leq i\leq 4$ and $$\psi_i=\left (t_i, \ \theta_{i}t_i + a_{i2}t_i^2+
\displaystyle\sum_{k\geq 3\atop i\not\in L_k}a_{ik}t_i^k\right )\ \ \ \mbox{for}\ \ \ 5\leq i\leq r.$$ The result follows since that $i\not\in L_k$ is equivalent to $k\leq \left [ \frac{i-2}{2}\right  ]$ for $k\geq 3$,.
\cqd

Notice that for $r\leq 3$ the moduli space corresponding is a single point and for $r=4$ is a one dimensional space.

For $r\geq 5$  the parameter vector associated to the normal form $\psi$ presented in the previous result is $\left (a_{5,2},a_{6,2},a_{6,3},\ldots ,a_{r,2},\ldots ,a_{r,\left [ \frac{r-2}{2}\right  ]}\right )$. Generically, $a_{5,2}\neq 0$ and by Theorem \ref{Anormal-form} we get a multigerm $\mathcal{A}$-equivalent to $\psi$ with parameter vector $p=\left (1,a'_{6,2},a'_{6,3},\ldots ,a'_{r,2},\ldots ,a'_{r,\left [ \frac{r-2}{2}\right  ]}\right )$.

In particular, the dimension of the moduli space $\mathbb{M}_r$ for multigerm of curves with singular ordinary point and $r\geq 5$ branches is the number of parameters in $p$ plus $r-3$ which corresponds to $\sharp\{\theta_i;\ 4\leq i\leq r\}$. If we denote $L_1=\{1,2,3\}$ just to unify the notation, we have that $\sharp I\setminus L_1 = \sharp\{\theta_i;\ 4\leq i\leq r\}$. Thus
$$\dim_{\mathbb{C}}\mathbb{M}_r=\sum_{k\geq 1}(\sharp I\setminus L_k)-1=\sum_{k=1}^{\left [ \frac{r-2}{2}\right  ]}(r-2k-1)=\left \{ \begin{array}{ll} \frac{(r-2)^2}{4} & \mbox{for}\ r\ even\\
\frac{(r-1)(r-3)}{4} & \mbox{for}\ r\ odd\end{array}\right .$$
as originally obtained by Granger in \cite{Granger}.

\subsection{Multigerms with $\Gamma_i\cap\mathbb{N}=\langle n,m\rangle$ and $\nu_i(f_l)=nm$}

We have just studied a class of multigerms in which
each block has just one element. In this subsection, we consider the opposite situation, that is, curves such that we have only one block.

Let $\mathcal{C}$ be a plane curve  defined by a reduced element
$f=\prod_{i=1}^{r}f_i\in\mathbb{C}\{X,Y\}$ such that each branch
$\mathcal{C}_i$ admits $\Gamma_i\cap\mathbb{N}=\langle n,m\rangle$,
$GCD(n,m)=1<n<m$ and $\nu_i(f_l)=nm$ for all $i,l\in I$ with $i\neq
l$. We will present a {\it pre-normal form} considering just the topological data and we will express the dimension of the generic component $\mathbb{M}_r(n,m)$ of the moduli space by means of the sets $L_k$. For the specific case of $n=2$ and $m=3$ we present a closed formula for $\dim_{\mathbb{C}}\mathbb{M}_r(2,3)$.

The curve $\mathcal{C}$ admits a multigerm with Puiseux block
form $\varphi=[\varphi_1,\ldots ,\varphi_r]$ where
$$\varphi_i=(t_i^n,a_{im}t_i^m+\sum_{j>m}a_{ij}t_i^j)$$ with
$a^n_{im}\neq a^n_{lm}, a_{im}\neq 0\neq a_{lm}$ for every $1\leq i,l\leq r$ and $i\neq l$. In this situation we have a single block and consequently the corresponding group to consider in Theorem \ref{A1normal-form} is $\mathcal{G}=\tilde{\mathcal{A}}_1$.

The topological class of $\mathcal{C}$ can be characterized by the value semiring $\Gamma$ that, according to Proposition 3.9 of \cite{carvalho1}, is $\Gamma=\langle v_1,v_2,\ldots ,v_{r+2}\rangle$ with
\begin{equation}\label{vi}
v_1=(n,\ldots ,n),\ \ v_2=(m,\ldots ,m),\ \ v_{i+2}=(nm,\ldots ,\infty,\ldots ,nm)\ \ i\in I,
\end{equation}
where the symbol $``\infty"$ occupies the $i^{th}$-coordinate in $v_{i+2}$.

The conductor of $\Gamma$ is $\kappa=(\kappa_1,\ldots ,\kappa_r)$ with $\kappa_i=rnm-n-m+1$ for all $i\in I$. In particular, if $k\geq rnm$ then $F^{\Gamma}_{J}(\underline{k})\neq \emptyset$ for every $\emptyset\neq J\subseteq I$ and as $\Gamma\subset\Gamma_1\times\ldots \times\Gamma_r$, if $k\not\in\langle n,m\rangle$, then $F^{\Gamma}_{J}(\underline{k})=\emptyset$.

The next lemma characterizes the non-empty fibers of $\underline{k}\in\mathbb{N}^r$ with respect to $\Gamma$ for $k<rnm$.

\begin{lemma}
    Given $k\in\langle n,m\rangle$ with $cnm\leq k<(c+1)nm$ and $0\leq c\leq r-1$ we have
    $$F^{\Gamma}_{J}(\underline{k})\neq\emptyset \Leftrightarrow \left \{ \begin{array}{lll} \sharp J\geq r-c & \mbox{if} & k-cnm\in\langle n,m\rangle \\ \sharp J\geq r-(c-1) & \mbox{if} & k-cnm\not\in\langle n,m\rangle. \end{array}\right .$$
\end{lemma}
\Dem Remark that $F^{\Gamma}_{J}(\underline{k})\neq\emptyset$ if and only if there exists $\gamma=(\gamma_1,\ldots ,\gamma_r)\in\Gamma$ with $\gamma_j=k$ for every $j\in J$ and $\gamma_j>k$ for every $j\in I\setminus J$. By (\ref{semiring}), $\gamma\in \Gamma$ can be expressed as
\begin{equation}\label{gamma}\gamma=\inf\left \{\sum_{s=1}^{r+2}\alpha_{1,s}v_s,\ldots ,\sum_{s=1}^{r+2}\alpha_{r,s}v_s \right \}\end{equation}
where $v_s=(v_{s1},\ldots ,v_{sr})$ is given in (\ref{vi}) for $1\leq s\leq r+2$.

For $j\in J$ there exists $l\in I$ such that $\gamma_j=k=\sum_{s=1}^{r+2}\alpha_{l,s}v_{sj}=\alpha_{l,1}n+\alpha_{l,2}m+\sum_{s=3}^{r+2}\alpha_{l,s}nm$ with $\alpha_{l,j+2}=0$.
In particular, $\sharp J\geq \sharp\{i\in I;\ \alpha_{l,i+2}=0\}=r-\sharp\{i\in I;\ \alpha_{l,i+2}\neq 0\}.$

If $k-cnm\not\in\langle n,m\rangle$ then
$\sum_{s=3}^{r+2}\alpha_{l,s}\leq c-1$ and $\sharp J\geq r-(c-1)$. If $k-cnm\in\langle n,m\rangle$, then $\sum_{s=3}^{r+2}\alpha_{l,s}\leq c$ and $\sharp J\geq r-c$.

On the other hand, if $k-cnm=\alpha_1n+\alpha_2m\in\langle
n,m\rangle$, considering any non negative integers $\alpha_3,\ldots
,\alpha_{r+2}$ such that $c=\sum_{s=3}^{r+2}\alpha_s$ we get
$\sharp\{i\in J;\ \alpha_{i+2}\neq 0\}\leq c$,
$\sum_{s=1}^{r+2}\alpha_sv_s\in F_{J}^{\Gamma}(\underline{k})$ for
$J=\{i\in I;\ \alpha_{i+2}=0\}$ and $\sharp J\geq r-c$. If
$k-cnm=g\not\in\langle n,m\rangle$, then
$k=g+nm+(c-1)nm=\alpha_1n+\alpha_2m+(c-1)nm$ and considering
$c-1=\sum_{s=3}^{r+2}\alpha_s$ we have
$\sum_{s=1}^{r+2}\alpha_sv_s\in F_{J}^{\Gamma}(\underline{k})$ for
$J=\{i\in I;\ \alpha_{i+2}=0\}$ and $\sharp J\geq r-(c-1)$. \cqd

Given $k>m$ we have that $\underline{k}\in\nu (\mathcal{M}^2)\subset\Lambda_{\tilde{\mathcal{A}_1}}$. So, considering $F_J(\underline{k})$ the $J$-fiber of $\underline{k}$ with respect to $\Lambda_{\tilde{\mathcal{A}_1}}$ we get $F^{\Gamma}_{J}(\underline{k+n})-\underline{n}\subset F_{J}(\underline{k})$ and, by the previous lemma, we obtain
    $$F^{\Gamma}_{J}(\underline{k+n})-\underline{n}\neq \emptyset\Leftrightarrow \left \{ \begin{array}{lll} \sharp J\geq r-c & \mbox{if} & k-cnm\in\langle n,m\rangle -n\\ \sharp J\geq r-(c-1) & \mbox{if} & k-cnm\not\in\langle n,m\rangle -n. \end{array}\right .$$

For each $k>m$ with $cnm-n\leq k<(c+1)nm-n$ taking
\begin{center}
$L=\{1,\ldots ,c+1\},\ J_i=\{i,c+2, c+3,\ldots ,r\}\ \ \mbox{for}\ \
i\in L\ \ \mbox{if}\ \ k-cnm\in\langle n,m\rangle -n;$

$L=\{1,\ldots ,c\},\ J_i=\{i,c+1, c+2,\ldots ,r\}\ \ \mbox{for}\ \
i\in L\ \  \mbox{if}\ \ k-cnm\not\in\langle n,m\rangle
-n$\end{center} in Proposition \ref{PropElim} and Corollary
\ref{CorElim} we obtain the following result:

\begin{proposition}\label{pre-normal}
    Any curve with value semiring generated by (\ref{vi}) is analytically equivalent to a curve defined by a multigerm $\varphi =[\varphi_1,\ldots ,\varphi_r]$ with
    $$\varphi_1=\left (t_1^n,a_{1m}t_1^m+\sum_{j>m\atop j\not\in\langle n,m\rangle-n}^{nm-n-1}a_{1j}t_1^j\right ),\ \ \ \ \varphi_i=\left (t_i^n,\sum_{j= m}^{(i-1)nm-n-1}a_{ij}t_i^j+\sum_{j\geq (i-1)nm-n\atop j-(i-1)nm\not\in\langle n,m\rangle-n}^{inm-n-1}a_{ij}t_i^j\right ),$$
    for $2\leq i\leq r$.
\end{proposition}

As in this topological class the value semiring $\Gamma$ is totally
determined by the $m$-jet of the multigerm $\varphi$, we can proceed
the study of $\Lambda_{\tilde{\mathcal{A}_1}}$ for the generic case
and consequently we get information about the generic component of
the moduli space $\mathbb{M}_r(n,m)$ of all plane curves with
$\Gamma$ generated by $(\ref{vi})$ considering the multigerm
$\varphi =[\varphi_1,\ldots ,\varphi_r]$ with
$\varphi_i=(t_i^n,\sum_{j\geq m}a_{ij}t_i^j)$ where
$a_{ij}\in\mathbb{C}$ are generic, $0\neq a_{im}^n\neq a_{lm}^n\neq
0$ for every $i,l\in I$ and $i\neq l$. The generic component for the
moduli space in this topological class was also considered in
\cite{genzmerpaul} by other methods.

Similarly to the case presented in the last subsection, the generality hypothesis implies that $L_k=\emptyset$, i.e., $e(k)=0$ or $L_k=\{1,\ldots ,e(k)\}$ with $e(k)\leq r$ for every $k> m$ and, by Remark \ref{lk-diff}, $e(k)$ is the maximum integer such that there exists $\rm{w}\in \mathcal{I}_{\mathcal{G}}$ with $j^k({\rm w})=(b_{1k}t_1^k,\ldots ,b_{rk}t_r^{k})$ for any $b_{ik}\in\mathbb{C}$ and $1\leq i\leq e(k)$.

In this context we have the following proposition:

\begin{proposition}\label{generic-dim}
The dimension of the generic component of the moduli space $\mathbb{M}_r(n,m)$ is zero for $r=1, (n,m)\in\{(2,m),(3,4),(3,5)\}$ and
$$\dim_{\mathbb{C}}\mathbb{M}_r(n,m)=r-2+\sum_{k> m}(r-e(k))$$
for the other cases.
\end{proposition}
\Dem Zariski, in \cite{zariskibook} proved that for irreducible
plane curves in the topological class determined by
$\Gamma\cap\mathbb{N}=\langle n,m\rangle$ with
$(n,m)\in\{(2,m),(3,4),(3,5)\}$ the moduli space is a single point.

By Theorem \ref{Anormal-form}, we can normalize the coefficient
$a_{1m}$ in $\varphi_1$ and thus the number of parameters in the
$\mathcal{A}$-normal form of $\varphi$ is $r-1+\sum_{k>m}\sharp
I\setminus L_k=r-1+\sum_{k> m}(r-e(k))$. Moreover, if $\sum_{k>
m}(r-e(k))\geq 1$ there exist $j=\min\{k>m;\ L_k\neq I\}$ and
$l=\min\{i\in I\setminus L_{j}\}$ such that $a_{lj}\neq 0$. In this
way, we can normalize such coefficient taking $c\in\mathbb{C}^*$
with $c^{j-m}=a_{lj}$, $\rho_i(t_i)=ct_i$ for all $i\in I$ and
$\sigma(X,Y)=(c^nX,c^mY)$, that is, $(\rho_1,\ldots
,\rho_r,\sigma)\cdot \varphi=\psi=[\psi_1,\ldots ,\psi_r]$ where the
coefficient of $t_{l}^{j}$ of $\psi_{l}$ is equal to $1$. In this
case, we obtain
$$\dim_{\mathbb{C}}\mathbb{M}_r(n,m)=r-2+\sum_{k> m}(r-e(k)).$$

We will show that, except for the above particular cases
considered by Zariski, we have $\sum_{k> m}(r-e(k))\geq 1$. More
specifically, we will show that $(r-e(m+1))+(r-e(m+2))\geq 1$.

Let us compute $e(m+1)$ that, as mentioned above, it is equivalent
to evaluate the maximum number of equations $j^{m+1}\left (
\frac{\varphi^*_i(\omega)}{nt_i^{n}}\right )=b_{i,m+1}t_{i}^{m+1}$
that admit solution for any $b_{i,m+1}\in\mathbb{C}$ with $0\leq
i\leq e(m+1)$ and $\omega=pdX+qdY\in\Omega_{\tilde{\mathcal{A}_1}}$,
that is, $p\in\langle X,Y\rangle^2$ and $q\in\langle X^2, Y\rangle$.

Notice that it is sufficient to consider $\nu_i\left
(\frac{\varphi^*_i(\omega)}{nt_i^n}\right )= m+1$. So, if $m\neq
n+1$ then $q=0$ and, if $m+1\not\equiv 0 \mod n$ then $p=0$. In
particular, for $m\neq n+1$ and $m+1\not\equiv 0 \mod n$ we have
$e(m+1)=0$.

For $m=n+1$ and $m+1\not\equiv 0\ mod\ n$ we get $\omega=\alpha
YdY$, $j^{m+1}\left ( \frac{\varphi^*_i(\omega)}{nt_i^{n}}\right
)=\frac{m}{n}\alpha a_{i,m}^2t_i^{m+1}$ and $e(m+1)=1$. Similarly,
if $m\neq n+1$ and $m+1\equiv 0\ mod\ n$ then $\omega
=\sum_{l=2}^{\frac{m+1}{n}}\alpha_lX^ldX$, $j^{m+1}\left (
\frac{\varphi^*_i(\omega)}{nt_i^{n}}\right
)=\sum_{l=2}^{\frac{m+1}{n}}\alpha_lt_i^{ln}$ and $e(m+1)=1$.

If $m=n+1$ and $m+1\equiv 0\ mod\ n$ then $n=2, m=3$ and we have
$\omega =\alpha_1 X^2dX+\alpha_2YdY$, that is, $j^4\left
(\frac{\varphi^*_i(\omega)}{2t_i^{2}}\right )=\left (
\alpha_1+\frac{3}{2}a_{i,3}^2\alpha_2\right )t_i^4$. So, we get
$e(m+1)=\min\{2,r\}$.

Now, we will compute $e(m+2)$ for the cases $e(m+1)>0$. As before it
is sufficient to consider
$\omega=pdX+qdY\in\Omega_{\tilde{\mathcal{A}_1}}$ with $\nu_i(p)\leq
m+2$ and $\nu_i(q)\leq n+2$.

In the case $m=n+1$ and $m+1\not\equiv 0\ mod\ n$, if $n>3$ we have
$\omega=\alpha YdY$, $j^{m+2}\left (
\frac{\varphi^*_i(\omega)}{nt_i^{n}}\right )=\frac{m}{n}\alpha
a^2_{i,m}t_i^{m+1}+\left ( \frac{2m+1}{n}\right )\alpha
a_{i,m}a_{i,m+1}t_i^{m+2}$
 and $e(m+2)=0$. If $n=3$ and thus $m=4$, we
have $\omega=\alpha_1 X^2dX+\alpha_2YdY$, $j^{6}\left (
\frac{\varphi^*_i(\omega)}{3t_i^{3}}\right )=\frac{4}{3}\alpha_2
a^2_{i,4}t_i^{5}+\left ( \alpha_1+3\alpha_2a_{i,4}a_{i,5}\right
)t_i^{6}$ and $e(m+2)=1$.

When $m\neq n+1$ and $m+1\equiv 0\ mod\ n$ we have the
possibilities: If $m\neq n+2$ and $n\neq 2$ then $\omega
=\sum_{l=2}^{\frac{m+1}{n}}\alpha_lX^ldX$, $j^{m+2}\left (
\frac{\varphi^*_i(\omega)}{nt_i^{n}}\right
)=\sum_{l=2}^{\frac{m+1}{n}}\alpha_lt_i^{ln}$, so $e(m+2)=0$. If
$m\neq n+2$ and $n=2$ then $\omega =\left
(\sum_{l=2}^{\frac{m+1}{2}}\alpha_lX^l+\alpha_1 XY\right
)dX+\alpha_0X^2dY$, $j^{m+2}\left (
\frac{\varphi^*_i(\omega)}{2t_i^{2}}\right
)=\sum_{l=2}^{\frac{m+1}{2}}\alpha_lt_i^{2l}+\left (
\alpha_1a_{i,m}+\alpha_0\frac{m}{2}a_{i,m}\right )t_i^{m+2}$ and
$e(m+2)=1$. If $m=n+2$ then $n=3, m=5$, $\omega=\alpha_1
X^2dX+\alpha_2YdY$ and $j^{7}\left (
\frac{\varphi^*_i(\omega)}{3t_i^{3}}\right
)=\alpha_1t^6_i+\frac{5}{3}\alpha_2a_{i,5}^2t_i^{7}$, consequently
$e(m+2)=1$.

Finally, if $n=2$ and $m=3$ then $\omega=(\alpha_1X^2+\alpha_2
XY)dX+(\alpha_3 Y+\alpha_4 X^2)dY$, $j^{5}\left (
\frac{\varphi^*_i(\omega)}{2t_i^{2}}\right )=\left
(\alpha_1+\frac{3}{2}\alpha_3 a^2_{i,3}\right )t^4_i+\left
(\alpha_2+\frac{3}{2}\alpha_4+\frac{7}{2}a_{i,4}\alpha_3\right
)a_{i,3}t_i^{5}$
 and $e(m+2)=1$.

These computations give us that, except to the cases $r=1$ and
$(n,m)\in\{(2,m),(3,4),(3,5)\}$, we have that
$(r-e(m+1))+(r-e(m+2))\geq 1$. \cqd

To illustrate the above result in a specific case, in what follows,
$n=2, m=3$ and $r\geq 2$.

We will compute $e(k)$ for any $k>m=3$ which allows us to exhibit
normal forms to the generic case and an explicit formula for
$\dim_{\mathbb{C}}\mathbb{M}_r(2,3)$.

Our strategy is the same as that of the previous proposition, that is,
to evaluate the maximum number of equations $j^{k}\left (
\frac{\varphi^*_i(\omega)}{2t_i^{2}}\right )=b_{i,k}t^{k}$ that
admit solution for any $b_{i,k}\in\mathbb{C}$ with
$\omega=pdX+qdY\in\Omega_{\tilde{\mathcal{A}_1}}$ that is,
$p\in\langle X,Y\rangle^2$ and $q\in\langle X^2, Y\rangle$. To
achieve this goal, for every $h\geq 2$ we consider
$$\omega_{2h}=\left ( \sum_{l=0}^{\left [ \frac{h}{3}\right ]}\alpha_{h-3l,2l}X^{h-3l}Y^{2l}\right )dX+\left ( \sum_{l=0}^{\left [ \frac{h-2}{3}\right ]}\beta_{h-2-3l,2l+1}X^{h-2-3l}Y^{2l+1}\right )dY$$
$$\omega_{2h+1}=\left ( \sum_{l=0}^{\left [ \frac{h-1}{3}\right ]}\alpha_{h-1-3l,2l+1}X^{h-1-3l}Y^{2l+1}\right )dX+\left ( \sum_{l=0}^{\left [ \frac{h}{3}\right ]}\beta_{h-3l,2l}X^{h-3l}Y^{2l}\right )dY.$$
Any $\omega\in\Omega_{\tilde{\mathcal{A}_1}}$ can be uniquely
expressed as $\omega=\sum_{l\geq 4}\omega_l$. As
\begin{center} $\nu_i\left (\frac{\varphi_i^*(X^{h-3l}Y^{2l}dX)}{2t_i^2} \right )=\nu_i\left (\frac{\varphi_i^*(X^{h-2-3l}Y^{2l+1}dY)}{2t_i^2} \right )=2h$ and
\vspace{0.4cm}

 $\nu_i\left (\frac{\varphi_i^*(X^{h-1-3l}Y^{2l+1}dX)}{2t_i^2} \right )=\nu_i\left (\frac{\varphi_i^*(X^{h-3l}Y^{2l}dY)}{2t_i^2} \right )=2h+1$
\end{center}
 for any $i\in I$, in order to evaluate $e(k)$ it is sufficient to consider $\omega=\sum_{l=4}^{k}\omega_l$.

\begin{proposition}\label{ek}
    If $\varphi=[\varphi_1,\ldots ,\varphi_r]$ is a generic multigerm
    admitting value semiring $\Gamma$ with generators as in (\ref{vi}), $n=2$, $m=3$ and $r\geq 2$ then $e(4)=2$, $e(5)=1$ and for $k\geq 6$:
    $$e(k)=\min\left \{2\left [\frac{k}{6}\right ]+1,r\right \}\ \mbox{if}\ k\not\equiv 4 \ mod\ 6, \ \ \ e(k)=\min\left \{2\left [\frac{k}{6}\right ]+3,r\right \}\ \mbox{if}\ k\equiv 4\ mod\ 6,$$
    and $e\left ( 6\left [ \frac{r-2}{2}\right ]+5\right )=e(3r-1)=r$ if $r$ is even.
\end{proposition}
\Dem  Remark that, by the last proposition, $e(4)=2$ and
$e(5)=1$.

Firstly we will consider the system $j^{c+1}\left
(\frac{\varphi^*_i(\omega_c+\omega_{c+1})}{2t_i^2} \right )=0$ for
every $i\in I$ and $c\geq 5$.

Using $\varphi_i=(t_i^2,\sum_{j\geq 3}a_{ij}t_i^j)$ we obtain
$$
\begin{array}{l}
j^{2h+1}\left ( \frac{\varphi^*_i(\omega_{2h})}{2t_i^2}\right )=\left ( \sum_{l=0}^{\left [ \frac{h}{3}\right ]}a_{i3}^{2l}\alpha_{h-3l,2l}+
\sum_{l=0}^{\left [ \frac{h-2}{3}\right ]}\frac{3}{2}a_{i3}^{2l+2}\beta_{h-2-3l,2l+1}\right )t_i^{2h}+\\
\\ \hspace{3cm} +\left ( \sum_{l=1}^{\left [ \frac{h}{3}\right
]}2la_{i3}^{2l-1}a_{i4}\alpha_{h-3l,2l}+ \sum_{l=0}^{\left [
\frac{h-2}{3}\right
]}\frac{6l+7}{2}a_{i3}^{2l+1}a_{i4}\beta_{h-2-3l,2l+1}\right
)t_i^{2h+1};
\end{array}$$
$$
\begin{array}{l}
j^{2h+2}\left ( \frac{\varphi^*_i(\omega_{2h+1})}{2t_i^2}\right )=\left ( \sum_{l=0}^{\left [ \frac{h-1}{3}\right ]}a_{i3}^{2l+1}\alpha_{h-1-3l,2l+1}+
\sum_{l=0}^{\left [ \frac{h}{3}\right ]}\frac{3}{2}a_{i3}^{2l+1}\beta_{h-3l,2l}\right )t_i^{2h+1}+\\
\\ \hspace{3cm} +\left ( \sum_{l=0}^{\left [ \frac{h-1}{3}\right
]}(2l+1)a_{i3}^{2l}a_{i4}\alpha_{h-1-3l,2l+1}+ \sum_{l=0}^{\left [
\frac{h}{3}\right ]}(3l+2)a_{i3}^{2l}a_{i4}\beta_{h-3l,2l}\right
)t_i^{2h+2}.
\end{array}$$

If $c=2h$ then to vanish the coefficient of $t_i^{c}$ in
$j^{c+1}\left (\frac{\varphi^*_i(\omega_{c}+\omega_{c+1})}{2t_i^2}
\right )$ for every $i\in I$ we obtain the system $N\cdot W_p^t=0$
with
$$ N=\left(\begin{array}{ccccc}
1 & a_{13}^2 & a_{13}^4 & \cdots & a_{13}^{2\left [ \frac{h}{3}\right ]+p}  \\
\vdots & \vdots & \vdots & & \vdots  \\
1 & a_{r3}^2 & a_{r3}^4 & \cdots & a_{r3}^{2\left [ \frac{h}{3}\right ]+p}  \\
\end{array}\right)\ \ \mbox{with}\ \ \left \{ \begin{array}{l} p=0,\ \ \mbox{if}\ c\equiv0, 2\ mod\
6;\\ p=2, \ \ \mbox{if}\  c\equiv4\ mod\ 6,\end{array}\right .$$ $$W_0=\left
(\alpha_{h,0}, \alpha_{h-3,2}+\frac{3}{2}\beta_{h-2,1},\ldots ,
\alpha_{h-3\left [ \frac{h}{3}\right ],2\left [ \frac{h}{3}\right
]}+\frac{3}{2}\beta_{h-3\left [ \frac{h}{3}\right ]+1,2\left [
\frac{h}{3}\right ]-1}\right )$$ and $W_2=\left
(W_0,\frac{3}{2}\beta_{h-2-3\left [ \frac{h}{3}\right ],2\left [
\frac{h}{3}\right ]+1}\right )$, where $W_p^t$ denotes the transpose
of $W_p$.

Considering the solution of $N\cdot W_p^t=0$, the associated system to vanish the coefficient of $t_i^{c+1}$ in $j^{c+1}\left ( \frac{\varphi_i^*(\omega_c+\omega_{c+1})}{2t_i^2}\right )$ for every $i\in I$ is $M\cdot Z^t = 0$ where
\begin{equation}\label{M}
M=\left(\begin{array}{cccccccc}
a_{14}a_{13} & a_{14}a_{13}^3 & \cdots & a_{14}a_{13}^{2\left [ \frac{h}{3}\right ]-1} & a_{13} & a_{13}^3 & \cdots & a_{13}^{2\left [ \frac{h}{3}\right ]+1}  \\
\vdots & \vdots &  & \vdots & \vdots & \vdots &  & \vdots \\
a_{r4}a_{r3} & a_{r4}a_{r3}^3 & \cdots & a_{r4}a_{r3}^{2\left [ \frac{h}{3}\right ]-1} & a_{r3} & a_{r3}^3 & \cdots & a_{r3}^{2\left [ \frac{h}{3}\right ]+1}  \\
\end{array}\right),\end{equation}
$$Z=\left(\frac{1}{2}\beta_{h-3l+1,2l-1};\ 1\leq l\leq \left [ \frac{h}{3}\right ], \alpha_{h-1-3l,2l+1}+\frac{3}{2}\beta_{h-3l,2l};\ 0\leq l\leq \left [ \frac{h}{3}\right ]-1,
\frac{3}{2}\beta_{h-3\left [ \frac{h}{3}\right ],2\left [ \frac{h}{3}\right ]}
\right)$$ if $c\equiv0\ mod\ 6$ or $Z=\left(\frac{1}{2}\beta_{h-3l+1,2l-1};\ 1\leq l\leq \left [ \frac{h}{3}\right ], \alpha_{h-1-3l,2l+1}+\frac{3}{2}\beta_{h-3l,2l};\ 0\leq l\leq \left [ \frac{h}{3}\right ]
\right)$ if $c\equiv2, 4\ mod\ 6$.

Notice that the rank of $M$ is $\min\{2\left [ \frac{h}{3}\right ]+1,r\}=\min\{2\left [ \frac{c}{6}\right ]+1,r\}$. If $2\left [ \frac{c}{6}\right ]+1\leq r$ then the solution of the above system implies that $\omega_{c}=0$.

For $c=2h+1$ we vanish the coefficient
of $t_i^{c}$ in $j^{c+1}\left (\frac{\varphi^*_i(\omega_{c}+\omega_{c+1})}{2t_i^2} \right )$ for every $i\in I$ solving the system $N\cdot W^t=0$ where
$$
N=\left(\begin{array}{cccc}
a_{13} & a_{13}^3 & \cdots & a_{13}^{2\left [ \frac{h}{3}\right ]+1}  \\
\vdots & \vdots &  & \vdots  \\
a_{r3} & a_{r3}^3 & \cdots & a_{r3}^{2\left [ \frac{h}{3}\right ]+1}  \\
\end{array}\right ),$$
\begin{center}
$W=(\alpha_{h-1-3l,2l+1}+\frac{3}{2}\beta_{h-3l,2l};\ 0\leq l\leq \left [ \frac{h}{3}\right ])$ if $c\equiv3, 5\ mod\ 6$ or

$W=(\alpha_{h-1-3l,2l+1}+\frac{3}{2}\beta_{h-3l,2l};\ 0\leq l\leq \left [ \frac{h}{3}\right ]-1,\frac{3}{2}\beta_{h-3\left [ \frac{h}{3}\right ],2\left [ \frac{h}{3}\right ]})$ if $c\equiv1\ mod\ 6$.
\end{center}

Taking the solution of the system $N\cdot W^t=0$ in order to vanish the coefficient of $t_i^{c+1}$ in $j^{c+1}\left ( \frac{\varphi^*_i(\omega_c+\omega_{c+1})}{2t_i^2}\right )$ for every $i\in I$ we obtain the system $M\cdot Z^t = 0$ where
$$
M=\left(\begin{array}{cccccccc}
a_{14} & a_{14}a_{13}^2 & \cdots & a_{14}a_{13}^{2\left [ \frac{h}{3}\right ]-2} & 1 & a_{13}^2 & \cdots & a_{13}^{2\left [ \frac{h}{3}\right ]}  \\
\vdots & \vdots &  & \vdots & \vdots & \vdots &  & \vdots \\
a_{r4} & a_{r4}a_{r3}^2 & \cdots & a_{r4}a_{r3}^{2\left [ \frac{h}{3}\right ]-2} & 1 & a_{r3}^2 & \cdots & a_{r3}^{2\left [ \frac{h}{3}\right ]}  \\
\end{array}\right)\ \mbox{if}\ c\equiv1\ mod\ 6, $$
$$M=\left(\begin{array}{cccccccc}
a_{14} & a_{14}a_{13}^2 & \cdots & a_{14}a_{13}^{2\left [ \frac{h}{3}\right ]} & 1 & a_{13}^2 & \cdots & a_{13}^{2\left [ \frac{h}{3}\right ]+2}  \\
\vdots & \vdots &  & \vdots & \vdots & \vdots &  & \vdots \\
a_{r4} & a_{r4}a_{r3}^2 & \cdots & a_{r4}a_{r3}^{2\left [ \frac{h}{3}\right ]} & 1 & a_{r3}^2 & \cdots & a_{r3}^{2\left [ \frac{h}{3}\right ]+2}  \\
\end{array}\right)\ \mbox{if}\ c\equiv3, 5\ mod\ 6,$$
$Z=\left(\frac{1}{2}\beta_{h-3l,2l};\ 0\leq l\leq \left [ \frac{h}{3}\right ]-1, \alpha_{h+1,0},\alpha_{h+1-3l,2l}+\frac{3}{2}\beta_{h+2-3l,2l-1};\ 1\leq l\leq \left [ \frac{h}{3}\right ]
\right)$ if $c\equiv1\ mod\ 6$,
$Z=\left(\frac{1}{2}\beta_{h-3l,2l};\ 0\leq l\leq \left [ \frac{h}{3}\right ], \alpha_{h+1,0},\alpha_{h+1-3l,2l}+\frac{3}{2}\beta_{h+2-3l,2l-1};\ 1\leq l\leq \left [ \frac{h}{3}\right ],\frac{3}{2}\beta_{h-1-3\left [ \frac{h}{3}\right ],2\left [ \frac{h}{3}\right ]+1}
\right)
$ if $c\equiv3\ mod\ 6$ and \newline
$Z=\left(\frac{1}{2}\beta_{h-3l,2l};\ 0\leq l\leq \left [ \frac{h}{3}\right ], \alpha_{h+1,0},\alpha_{h-2-3l,2l+2}+\frac{3}{2}\beta_{h-1-3l,2l+1};\ 0\leq l\leq \left [ \frac{h}{3}\right ]
\right)
$ if $c\equiv5\ mod\ 6$.

Remark that the rank of $M$ is $\min\{2\left [ \frac{h}{3}\right ]+1,r\}=\min\{2\left [ \frac{c}{6}\right ]+1,r\}$ if $c\equiv1\ mod\ 6$ or $\min\{2\left [ \frac{h}{3}\right ]+3,r\}=\min\{2\left [ \frac{c}{6}\right ]+3,r\}$ $c\equiv3, 5\ mod\ 6$. In particular $\omega_{c}=0$ if $2\left [ \frac{c}{6}\right ]+1\leq r$ for $c\equiv1\ mod\ 6$ or $2\left [ \frac{c}{6}\right ]+3\leq r$ for $c\equiv 3, 5\ mod\ 6$.

In this way to study the system $j^k\left (\frac{\varphi^*_i(\omega)}{2t_i^2}\right )=b_{i,k}t_i^k$ with $b_{i,k}\in\mathbb{C}$, $i\in I$ for $r\leq 2\left [ \frac{k-1}{6}\right ]+1$ or $2\left [ \frac{k-1}{6}\right ]+3\leq r$ it is sufficient to consider $\omega=\omega_{k-1}+\omega_k\in\Omega_{\tilde{\mathcal{A}_1}}$. For this situation the cases are the same as those considered above where we obtain the system $M\cdot Z^t=(b_{1,k},\ldots ,b_{r,k})^t$ with $M$ and $Z$ previously described taking $k=c+1$. Hence, $e(k)=rank(M)$ for all $k\geq 6$.

As $\left [ \frac{k}{6}\right ]=\left [ \frac{c}{6}\right ]$ for $c\not\equiv 5\ mod\ 6$ and $\left [ \frac{k}{6}\right ]=\left [ \frac{c}{6}\right ]+1$ for $c\equiv 5\ mod\ 6$ we get
$$e(k)=\min\left \{2\left [\frac{k}{6}\right ]+1,r\right \}\ \mbox{if}\ k\not\equiv 4 \ mod\ 6 \ \ \ \mbox{and}\ \ \ e(k)=\min\left \{2\left [\frac{k}{6}\right ]+3,r\right \}\ \mbox{if}\ k\equiv 4\ mod\ 6.$$

We obtain the same conclusion if $r=2\left [ \frac{k-1}{6}\right ]+2$ and $k\not\equiv 5\ mod\ 6$.

For $r=2\left [ \frac{k-1}{6}\right ]+2$ and $k\equiv 5\ mod\ 6$, that is, $k=3r-1=6\left [ \frac{r-2}{2}\right ]+5$, the condition $j^{k-1}\left ( \frac{\varphi_i^*(\omega_{k-2}+\omega_{k-1})}{2t_i^2}\right )=0$ for every $i\in I$ produces a system $M\cdot Z^t=0$ with $rank(M)=\min\{2\left [ \frac{k-2}{6}\right ]+3,r\}=\min\{r+1,r\}=r$. In this way, we have an extra variable for the system $j^{k}\left ( \frac{\varphi_i^*(\omega_{k-2}+\omega_{k-1}+\omega_k)}{2t_i^2}\right )=b_{i,k}t_i^k$ that can be expressed as $M_1\cdot Z_1^t=(b_{1,k},\ldots ,b_{r,k})^t$ where $M_1$ is the matrix $M$ as in (\ref{M}) with an extra column depending on $a_{i3}, a_{i4}$ and $a_{i5}$. This allows us to conclude that $$e\left ( 6\left [ \frac{r-2}{2}\right ]+5\right )=e(3r-1)=rank(M_1)=\min\left \{2\left [ \frac{k-1}{6}\right ]+1+1,r\right \}=r.$$
\cqd

As a consequence of Propositions \ref{generic-dim} and \ref{ek} we obtain an explicit  formula for the dimension of the generic component of $\mathbb{M}_r(2,3).$

\begin{corollary}\label{23} For plane curves that admit value semiring generated by (\ref{vi}) with $n=2$, $m=3$ and $r$ branches the dimension of the generic component $\mathbb{M}_r(2,3)$ of the moduli space is
$$\dim_{\mathbb{C}}\mathbb{M}_r(2,3)=\left \{\begin{array}{ll}
\frac{(r-1)(3r-5)}{2}& \mbox{if}\ r\ \mbox{is odd} \vspace{0.2cm}\\
\frac{(r-1)(3r-5)+1}{2}& \mbox{if}\ r\ \mbox{is even.}
\end{array}\right .$$
\end{corollary}
\Dem The case $r=1$ is immediate.

By Proposition \ref{generic-dim} and the above result we have that
\begin{equation}\label{dim}\begin{array}{l}\dim_{\mathbb{C}}\mathbb{M}_r(2,3)=r-2+\displaystyle{\sum_{k\geq 4}(r-e(k))}=3r-5+\displaystyle{\sum_{k\geq 6}(r-e(k))}=\\
\hspace{2.4cm} =3r-5+5\hspace{-0.3cm}\displaystyle{\sum_{k\geq 6\atop k\not\equiv 4\ mod\ 6}\hspace{-0.4cm}(r-e(k))}+\hspace{-0.4cm}\displaystyle{\sum_{k\geq 6\atop k\equiv 4\ mod\ 6}\hspace{-0.4cm}(r-e(k))}.\end{array}\end{equation}

Remark that for $r=2$ we have $e(k)=2$ for every $k\geq 6$ then $\dim_{\mathbb{C}}\mathbb{M}_2(2,3)=1$.

If $k\equiv 4\ mod\ 6$ we have $e(k)=\min\left \{2\left [\frac{k}{6}\right ]+3,r\right \}$ and \begin{equation}\label{k2}\displaystyle{\sum_{k\geq 6\atop k\equiv 4\ mod\ 6}\hspace{-0.3cm}(r-e(k))}=\sum_{i=1}^{\left [ \frac{r-3}{2}\right ]}(r-(2i+3))=\left [ \frac{r-3}{2}\right ]\left ( r-4- \left [ \frac{r-3}{2}\right ]\right ).\end{equation}

For $k\not\equiv 4\ mod\ 6$ then $e(k)=\min\left \{2\left [\frac{k}{6}\right ]+1,r\right \}$ and $e\left ( 6\left [ \frac{r-2}{2}\right ]+5\right )=e(3r-1)=r$ if $r$ is even. So,
{\small \begin{equation}\label{k1}\displaystyle{\sum_{k\geq 6\atop k\not\equiv 4\ mod\ 6}\hspace{-0.3cm}(r-e(k))}=\left \{\begin{array}{ll}
\sum_{i=1}^{\left [ \frac{r-4}{2}\right ]}(r-(2i+1))=\left [ \frac{r-1}{2}\right ]\left ( r-2- \left [ \frac{r-1}{2}\right ]\right )-1 & \mbox{if}\ k\equiv 5\ mod\ 6\ \mbox{and}\ r\ \mbox{even}\vspace{0.2cm}\\
\sum_{i=1}^{\left [ \frac{r-1}{2}\right ]}(r-(2i+1))=\left [ \frac{r-1}{2}\right ]\left ( r-2- \left [ \frac{r-1}{2}\right ]\right ) & \mbox{otherwise}.
\end{array}
\right .\end{equation}
}

For $r\geq 3$, considering (\ref{k2}) and (\ref{k1}) in (\ref{dim}) the result follows.
\cqd

{\bf Acknowledgment:}
We would like to thank Professor Yohann Genzmer for valuable discussions and we express our sincere gratitude to Professor Abramo Hefez and the anonymous referee for their careful reading of the manuscript and their comments and suggestions.

\vspace{0.5cm}

\begin{center}
\begin{tabular}{ccc}
Hernandes, M. E. & & Rodrigues Hernandes, M. E. \\
{\it mehernandes$@$uem.br} & & {\it merhernandes$@$uem.br}\\
\end{tabular}
\vspace{0.5cm}

Universidade Estadual de Maring\'a

Maring\'a - Paran\'a - Brazil
\end{center}

\end{document}